\documentclass[12pt,oneside,reqno]{amsart}

\usepackage[T1]{fontenc}
\usepackage[utf8]{inputenc}
\usepackage{amsmath, amssymb, amscd, amsthm, amsfonts}
\usepackage{newtxtext,newtxmath}

\usepackage{fancybox}
\usepackage{graphicx}
\usepackage{eepic}
\usepackage{color}
\usepackage{enumerate}
\usepackage{comment}
\usepackage{mathtools}
\usepackage{mathrsfs}
\usepackage[all,cmtip]{xy}
\usepackage{tikz-cd}
\usepackage{wrapfig}
\usepackage{enumitem}
\usepackage{stmaryrd}
\usepackage{tikz}
\usetikzlibrary{arrows.meta,positioning}
\usepackage{tkz-euclide}
\usepackage[hidelinks]{hyperref}

\newtheorem{theorem}{Theorem}[section]

\newtheorem{lemma}[theorem]{Lemma}

\newtheorem{prop}[theorem]{Proposition}

\newtheorem{defin}[theorem]{Definition}

\newtheorem{remark}[theorem]{Remark}

\newtheorem{corollary}[theorem]{Corollary}

\newcommand{\C}{\mathbb{C}}
\newcommand{\R}{\mathbb{R}}
\newcommand{\Q}{\mathbb{Q}}
\newcommand{\A}{\mathbb{A}}

\newcommand{\Z}{\mathbb{Z}}

\newcommand{\F}{\mathbb{F}}
\renewcommand{\O}{\mathcal{O}}

\newcommand{\Gal}{{\rm Gal}}

\newcommand{\plim}[1][]{\mathop{\varprojlim}\limits_{#1}}

\def\kakko<#1>{\langle #1 \rangle}
\def\hoge[#1]{\llbracket #1 \rrbracket}

\title{Explicit Galois Deformations over Imaginary Quadratic Fields}
\author{Takuya Tanaka}
\date{}
\address{Takuya Tanaka \\
Department of Mathematics, Institute of Science Tokyo\\
 2-12-1, Ookayama, Meguro, Tokyo, 152-8550, JAPAN}
\email{tanaka.t.ch@m.titech.ac.jp}
\begin{document}

\begin{abstract}
Let $K$ be an imaginary quadratic field, and let $p$ be an odd prime that splits as
$(p)=\pi\bar{\pi}$ in $K$. Let $G_K^\pi$ and $G_K^p$ denote the Galois groups of the maximal algebraic extensions of $K$ unramified outside $\pi$ and outside the primes above $p$, respectively. In this paper, we construct certain Iwasawa-theoretic quotients of 
$G_K^\pi$ and $G_K^p$ and
study the universal deformation rings of the induced residual representations.
The defining relations of these quotient deformation rings are described
through characteristic elements of Iwasawa modules, and hence through Katz
$p$-adic $L$-functions.

\end{abstract}
\maketitle

\tableofcontents

\section{Introduction}\label{section-introduction}
The study of universal deformation rings of mod $p$ Galois representations of number
fields, with prescribed ramification, is a central topic in number theory.  Let
$G_K^{S}:=\mathrm{Gal}(K^{S}/K)$ be the Galois group of the maximal algebraic
extension of a number field $K$ unramified outside a finite set $S$ of finite places.
Given a residual representation
$\bar\rho:G_K^{S}\to \mathrm{GL}_2(\mathbb F_p)$, one may ask how explicitly its
universal deformation ring can be described.

In this paper, we study this question when $K$ is an imaginary quadratic field and
$S=\{\pi\}$ or $S=\{\pi,\overline{\pi}\}$, where $p$ splits in $K$ as
$(p)=\pi\overline{\pi}$.  In the final section, we apply our results to Galois representations
arising from Bianchi cusp forms.  We construct a large Iwasawa-theoretic quotient
of $G_K^{S}$ through which $\bar\rho$ factors, and we describe explicitly both the
universal deformation ring of the induced representation of this quotient and the
universal deformation itself on generators.  The defining relations are expressed in
terms of characteristic elements of Iwasawa modules and hence, via Rubin's main
conjecture, in terms of Katz $p$-adic $L$-functions.  Our approach follows the work
of Boston \cite{BostonThesis} and M\'ezard \cite{Mezard99}, who treated the case $K=\Q$.
The main new feature in the present setting is that the relations in the deformation
ring are obtained from characteristic elements of Iwasawa modules over imaginary
quadratic fields.  Moreover, we will explicitly describe a universal deformation of a Galois representation attached to a Bianchi cusp form considered by
Berger--Klosin \cite{BergerKlosin}.

Throughout this paper, let $p$ be an odd prime which splits in $K/\Q$ as
$(p)=\pi\overline{\pi}$.  For $\bullet\in\{\pi,p\}$, put
\[
S_\pi:=\{\pi\},\qquad S_p:=\{\pi,\overline{\pi}\},
\]
and write $G_K^\bullet:=G_K^{S_\bullet}$.  Let
$\bar\rho:G_K^{S_\bullet}\to \mathrm{GL}_2(\F_p)$ be a Galois representation whose
image is contained in the subgroup of upper triangular matrices but is not diagonal.
We write
\[
\bar\rho=
\begin{pmatrix}
\chi_1 & * \\
0 & \chi_2
\end{pmatrix}
\]
for its diagonal characters $\chi_1$ and $\chi_2$.  Let $\mathcal C$ be the category
of complete Noetherian local rings with residue field $\F_p$.  Consider the functor
\[
\mathcal F_{\bar\rho}:\mathcal C\longrightarrow \mathrm{Set},\qquad
R\longmapsto
\left\{
\rho_R:G_K^{S_\bullet}\to \mathrm{GL}_2(R)\ \middle|\ 
\rho_R\bmod \mathfrak m_R=\bar\rho
\right\}/\sim .
\]
Here $\mathfrak m_R$ is the maximal ideal of $R$, and $\rho_1\sim\rho_2$ means that
$\rho_1=M\rho_2M^{-1}$ for some
$M\in \ker\{\mathrm{GL}_2(R)\to \mathrm{GL}_2(\F_p)\}$.

We assume that $\operatorname{End}_{\F_p[G_K^{S_\bullet}]}(\bar\rho)=\F_p$.  Under this assumption,
Mazur \cite{Mazur} proved that $\mathcal F_{\bar\rho}$ is representable.  Let
$\mathcal R_{\bar\rho}$ be a representing object.  Thus, for all $T\in\mathcal C$,
there is a canonical isomorphism
\[
\mathcal F_{\bar\rho}(T)\simeq
\operatorname{Hom}_{\mathcal C}(\mathcal R_{\bar\rho},T).
\]
We call $\mathcal R_{\bar\rho}$ the universal deformation ring of $\bar\rho$.

Let $L_{\bar\rho}:=\overline K^{\ker\bar\rho}$ be the fixed field of $\ker\bar\rho$, and let
$K\subset F\subset L_{\bar\rho}$ be the maximal abelian subextension of $K$.  Then $[F:K]$ is
prime to $p$.  Let $F^{S_\bullet}(p)$ be the maximal pro-$p$ extension of $F$
unramified outside the primes above $S_\bullet$.  Since
$\ker(\mathrm{GL}_2(R)\to \mathrm{GL}_2(R/\mathfrak m_R))$ is a pro-$p$ group for
all $R\in\mathcal C$, every deformation of $\bar\rho$ factors through
$\mathrm{Gal}(F^{S_\bullet}(p)/K)$.  In particular, the universal deformation induces
a representation
\[
\rho^{\mathrm{univ}}:\mathrm{Gal}(F^{S_\bullet}(p)/K)\to
\mathrm{GL}_2(\mathcal R_{\bar\rho}).
\]

Let $F_\infty/F$ be the maximal multiple $\Z_p$-extension unramified outside the
primes above $S_\bullet$ and abelian over $K$.  Let $M$ be the maximal abelian
subextension of $F^{S_\bullet}(p)/F_\infty$, and define the Iwasawa-theoretic quotient
\[
\mathrm{Gal}(F^{S_\bullet}(p)/K)\twoheadrightarrow \mathrm{Gal}(M/K)=:G^{\bullet}.
\]
We study the universal deformation of the induced representation
$\bar\rho:G^{\bullet}\to \mathrm{GL}_2(\F_p)$.  This type of quotient deformation
problem was considered by Boston and M\'ezard.

Before stating the main theorem, we illustrate it in a simple case.  In the following, we always use $W$ to denote the Witt ring $W=W(\bar\F_p)$. Let
\(\mu_F(\chi)\in W\hoge[T]\) and
\(\mu_{F,\pi}(\chi)\in W\hoge[S,T]\) denote, respectively, the one-variable and
two-variable Katz \(p\)-adic \(L\)-functions associated with the choice of
\(\pi\) and a character \(\chi\in\widehat{\Delta}\)
(Definition \ref{def:p-L}).  Let $F/K$ be
the abelian extension defined above, and let $\mathfrak{X}_{F,\pi}$ be the Iwasawa module unramified
outside $\pi$.  Assume that  $1\neq \chi_0:=\chi_1/\chi_2$ is non-quadratic and the $\chi_0$-component of $\mathfrak{X}_{F,\pi}$ is cyclic as a
$\Lambda$-module, that $\mathrm{Cl}_K\otimes\Z_p$ is nonzero and cyclic, and that the
maximal $\Z_p$-extension unramified outside $\pi$ contains no nontrivial unramified
subextension.  Then the quotient deformation ring $R_G$ has the following simple
presentation:
\[
R_G\simeq
\Z_p\hoge[S,T,V]
/
\left(
 f\!\left(\frac{1+S}{1+T}-1\right),\ (1+V)^a-1
\right),
\]
where $f\in\Z_p\hoge[T]$ satisfies $(f)=(\mu_F(\chi_0))$ in
$W\hoge[T]$, and $a=p^{\mathrm{val}_W(\mu_K(0))}$, where \(\mu_K(0)\) is the constant term of the corresponding Katz
\(p\)-adic \(L\)-function attached to \(\pi\) over \(K\) and $\mathrm{val}_W$ is the normalized valuation of $W$.  By universality, we obtain a natural
surjection $\mathcal R_{\bar\rho}\twoheadrightarrow R_G$.  We also study conditions
under which this surjection is an isomorphism (Corollary \ref{cor:isom-condition}). We shall distinguish integral coefficients from \(W\)-coefficients.
Put
\[
\Lambda:=\mathbb Z_p\hoge[\Gamma],\qquad
\Lambda_W:=\Lambda \widehat{\otimes}_{\mathbb Z_p}W,
\]
where $\Gamma:=\Gal(F_\infty/F)$. Whenever topological generators
$\gamma_1,\ldots,\gamma_d$ of $\Gamma$ are fixed, we identify $\Lambda$ with the
formal power series ring $\Z_p\hoge[T_1,\ldots,T_d]$ by
$\gamma_i-1\mapsto T_i$. For a matrix $M\in M_{a,b}(R)$ over a commutative
ring $R$, we denote by $I_k(M)$, $k\leq \min\{a,b\}$ the ideal of $R$ generated by the determinants
of all $k\times k$ submatrices of $M$. For the above $F/K$, put $P_{F,\bullet}:=\Gal(M/F)\subset G^{\bullet}$ for each $\bullet \in \{p,\pi\}$.
\begin{theorem}[Main Theorem]\label{Main 2}
Let $K$ be an imaginary quadratic field.  Let
$\bar\rho: G_K^\bullet\to \mathrm{GL}_2(\F_p)$, $\bullet\in\{\pi,p\}$, be a Galois
representation whose image is noncommutative and whose diagonal characters
$\chi_1,\chi_2$ satisfy $\chi_0:=\chi_1/\chi_2\neq 1$.  Assume that
$\mathrm{Cl}_K\otimes \Z_p$ is cyclic and that $\chi_0$ is not quadratic.

We construct an Iwasawa-theoretic quotient
$G_K^\bullet\twoheadrightarrow G^\bullet$
through which $\bar\rho$ factors.  Let
$\bar\rho_G:G^\bullet\to \mathrm{GL}_2(\F_p)$ be the induced representation.  Let
$R_G$ be its universal deformation ring, and let
$\rho_G^{\mathrm{univ}}:G^\bullet\to \mathrm{GL}_2(R_G)$ be the universal deformation.
Then $R_G$ has an explicit presentation, and $\rho_G^{\mathrm{univ}}$ is explicitly
described on generators as follows.

\begin{enumerate}
\item Suppose that $\bullet=\pi$.  The pro-$p$ group $P_{F,\pi}\subset G^\pi$ admits a minimal system of
generators
\[
\langle s_1,\ s_{\chi,1},\ldots,s_{\chi,s_\chi},\ \gamma
\mid 1\neq\chi\in\widehat{\Delta}\rangle
\]
such that $\rho_G^{\mathrm{univ}}(s_{\chi,i})=1$ if $\chi\neq 1,\chi_0$, and
\[
\rho_G^{\mathrm{univ}}(s_{\chi_0,i})=
\begin{pmatrix}
1 & U_i\\
0 & 1
\end{pmatrix}\ (i\ge 2),\quad
\rho_G^{\mathrm{univ}}(\gamma)=
\begin{pmatrix}
1+S & 0\\
0 & 1+T
\end{pmatrix},
\]
\[
\rho_G^{\mathrm{univ}}(s_{\chi_0,1})=
\begin{pmatrix}
1 & 1\\
0 & 1
\end{pmatrix},\quad
\rho_G^{\mathrm{univ}}(s_1)=
\begin{pmatrix}
1+V & 0\\
0 & 1+V
\end{pmatrix}.
\]
Let $a_\pi:=p^{\mathrm{val}_W(\mu_K(0))}$, let $\delta\in\{0,1\}$, where
$\delta=1$ if $\mathrm{Cl}_K\otimes \Z_p=0$ or $K_\infty$ contains the maximal unramified $p$-extension of $K$, and $\delta=0$ otherwise. There exists a matrix
$F_{\chi_0}=(f_{k\ell})_{k,\ell}\in M_{s_{\chi_0},s_{\chi_0}}(\Lambda)$ such that, after extension of scalars to $\Lambda_W$,
\[
\det(F_{\chi_0})\Lambda_W=(\mu_F(\chi_0)),
\]
and
\begin{align*}
R_G\simeq
\Z_p\hoge[S,T,U_2,\ldots,U_{s_{\chi_0}},V]
\bigg/
\bigg(
& f_{1,\ell}\!\left(\frac{1+S}{1+T}-1\right)
 +\sum_{k=2}^{s_{\chi_0}}U_k
 f_{k,\ell}\!\left(\frac{1+S}{1+T}-1\right),\\
& (1+V)^{a_\pi}-1,\ \delta V
\ \bigg|\
\ell=1,\ldots,s_{\chi_0}
\bigg).
\end{align*}

\item Suppose that $\bullet=p$.  The pro-$p$ group $P_{F,p}\subset G^p$ admits a minimal system of generators
\[
\langle \gamma_1,\gamma_2,\ s_1,\ s_{\chi,1},\ldots,s_{\chi, s_\chi}
\mid 1\neq\chi\in\widehat{\Delta}\rangle
\]
such that $\rho_G^{\mathrm{univ}}(s_{\chi,i})=1$ if $\chi\neq 1,\chi_0$, and
\[
\rho_G^{\mathrm{univ}}(s_{\chi_0,i})=
\begin{pmatrix}
1 & U_i\\
0 & 1
\end{pmatrix}\ (i\ge 2),\quad
\rho_G^{\mathrm{univ}}(\gamma_i)=
\begin{pmatrix}
1+S_i & 0\\
0 & 1+T_i
\end{pmatrix},
\]
\[
\rho_G^{\mathrm{univ}}(s_{ \chi_0,1})=
\begin{pmatrix}
1 & 1\\
0 & 1
\end{pmatrix},\quad
\rho_G^{\mathrm{univ}}(s_1)=
\begin{pmatrix}
1+V & 0\\
0 & 1+V
\end{pmatrix}.
\]
Let $\delta\in\{0,1\}$, where $\delta=1$ if
$\mathrm{Cl}_K\otimes \Z_p=0$ or $K_\infty$ contains the maximal unramified
$p$-extension of $K$, and $\delta=0$ otherwise.  There exists an element $a\in\Z_p$ and a presentation matrix
$G_{\chi_0}=(g_{k\ell})\in M_{s_{\chi_0},s_{\chi_0}-1}(\Lambda)$ such that, after extension of scalars to $\Lambda_W$, for some non-zero ideals $\mathfrak{a}, \mathfrak{a}^{'}\subset \Lambda_W$,
\[
I_{s_{\chi_0}-1}(G_{\chi_0})\Lambda_W
\supset
\mathfrak{a}\mu_{F,\pi}(\chi_0)+\mathfrak{a}^{'}\mu_{F,\bar\pi}(\chi_0),
\]
and
\begin{align*}
R_G\simeq
\Z_p&\hoge[S_1,T_1,S_2,T_2,U_2,\ldots,U_{s_{\chi_0}},V]
\bigg/
\bigg(
 g_{1,\ell}\!\left(\frac{1+S_1}{1+T_1}-1,
                    \frac{1+S_2}{1+T_2}-1\right)\\
&+\sum_{k=2}^{s_{\chi_0}}U_k
 g_{k,\ell}\!\left(\frac{1+S_1}{1+T_1}-1,
                    \frac{1+S_2}{1+T_2}-1\right),
\ (1+V)^a-1,\ \delta V
\ \bigg|\
\ell=1,\ldots,s_{\chi_0}-1
\bigg).
\end{align*}
\end{enumerate}
\end{theorem}
Although our main results concern Iwasawa-theoretic quotient deformation
rings, in certain special cases the same method determines the full universal deformation ring. This occurs, for instance, when the relevant pro-\(p\)
Galois group over \(F\) is free.
We also study the condition under which the ideals $\mathfrak{a}, \mathfrak{a}^{'}$ are equal to $(1)$ (Corollary \ref{cor:cyclicity}).
For $p=3,5,7$, the imaginary quadratic fields $\Q(\sqrt{-d})$, $d\leq 100$, for which
the $p$-part of the ideal class group is cyclic and nontrivial are as follows:
\[
\langle d,p\rangle=
\langle 23,3\rangle,\ \langle 26,3\rangle,\ \langle 29,3\rangle,\ \langle 31,3\rangle,
\langle 38,3\rangle,\ \langle 53,3\rangle,\ \langle 59,3\rangle,\ \langle 61,3\rangle,
\]
\[
\langle 83,3\rangle,\ \langle 87,3\rangle,\ \langle 89,3\rangle,
\langle 47,5\rangle,\ \langle 74,5\rangle,\ \langle 79,5\rangle,
\langle 86,5\rangle,\ \langle 71,7\rangle.
\]

We now outline the proof.  First, we construct an Iwasawa-theoretic quotient of
$G_K^\pi$ and study its topological generators and relations.  Second, we relate the
defining relations to characteristic elements of Iwasawa modules
(Theorem \ref{Main 1}).  By Rubin's main conjecture, these characteristic elements are
expressed in terms of Katz $p$-adic $L$-functions associated with \(\pi\).  Finally, we compute the images of
the generators and relations under the universal deformation.  The same strategy also
applies in the $p$-ramified case.

One motivation for this work comes from the deformation theory of mod $p$ Galois
representations attached to Bianchi cusp forms at Eisenstein primes.  In the final
Section, we treat such modular examples.  In particular, we treat the universal
deformation of the mod $5$ Galois representation over
$K=\Q(\sqrt{-51})$ considered by Berger--Klosin \cite{BergerKlosin}
(Proposition~\ref{prop:autom}).

The structure of this paper is as follows.
\begin{itemize}
\item In Section 2, we briefly review Greenberg's work \cite{Gr} on pseudo-null
submodules of finitely generated $\Lambda$-modules.
\item In Section 3, following the method of B\"olling \cite{Bo}, we investigate the
structure of the Galois group of the maximal pro-$p$ extension of an imaginary
quadratic field unramified outside $\pi$ and $p$.
\item In Section 4, we recall the Iwasawa
main conjecture for imaginary quadratic fields.  We then study the structure of
Iwasawa modules.
\item In Section 5, we relate the relation of the Galois group to the $p$-adic $L$-function associated with \(\pi\).
\item In Sections 6 and 7, we explain an application of our results to mod $p$ Galois
representations.
\end{itemize}

\section*{Acknowledgements}
I would like to express my gratitude to JSPS for its financial support.

\section{Finitely Generated $\Lambda$-Modules}
Let $k$ be a number field and let $S$ be a set of finite places of $k$. Put $k_\infty$ be the maximal multiple $\Z_p$-extension unramified outside $S$ and put $\mathfrak{X}_{k,S}$ be the Galois group of the maximal abelian pro-$p$ extension of $k_\infty$ unramified outside $S$ over $k_\infty$. 
In this section, we study pseudo-null submodules of $\mathfrak{X}_{k,S}$. Let $k_\infty/k$ be the field extension with
Galois group $G\simeq \mathbb Z_p^d$. 

\begin{prop}\label{lem:fg}
$\mathfrak{X}_{k,S}$ is a finitely generated $\Lambda_G$-module.
\end{prop}

The key result is the following topological Nakayama lemma.

\begin{lemma}\label{lem:top-nakayama}
Let $X$ be a compact $\Lambda_G$-module. The following are equivalent:
\begin{itemize}
    \item $X$ is finitely generated over $\Lambda_G$.
    \item $X/\mathfrak m_GX$ is finite, where $\mathfrak m_G$ is the maximal
    ideal of $\Lambda_G$.
\end{itemize}
\end{lemma}

\begin{proof}[Proof of Proposition~\ref{lem:fg}]
In the following proof, for each subfield $L$ of the $k_\infty/k$, we denote $\mathfrak{X}(L)$ be the Galois group of the maimal abelian pro-$p$ extension unramified outside $S$ for simplicity. So, $\mathfrak{X}(k_\infty)=\mathfrak{X}_{k,S}$.
We first treat the case $d=1$. Let $g$ be a topological generator of $G$.
By assumption, the image of the  natural inclusion
$\mathfrak{X}(k_\infty)/(g-1)\mathfrak{X}(k_\infty)\to \mathfrak{X}(k)$ is a finitely generated
$\Z_p$-module due to the class field theory. 
In particular, $\mathfrak{X}(k_\infty)/\mathfrak m_G\mathfrak{X}(k_\infty)$ is finite. Hence
$\mathfrak{X}(k_\infty)$ is finitely generated over $\Lambda_G$ by
Lemma \ref{lem:top-nakayama}.

Now assume $d>1$. Choose a direct summand
$H=\langle h\rangle\simeq\mathbb Z_p$ of $G$, and put
$G':=G/H$ and $k_\infty':=k_\infty^H$. Then
$\Lambda_G/(h-1)\Lambda_G\simeq \Lambda_{G'}$. By the case $d=1$ applied
relatively to $k_\infty/k_\infty'$, the module
$\mathfrak{X}(k_\infty)/(h-1)\mathfrak{X}(k_\infty)$ is finitely generated over $\Lambda_{G'}$.
By induction on $d$, this implies that $\mathfrak{X}(k_\infty)$ is finitely generated
over $\Lambda_G$.
\end{proof}

If $G=\langle\gamma_1,\ldots,\gamma_d\rangle\simeq \mathbb Z_p^d$, then
Serre's isomorphism gives
\[
\Lambda_G=\mathbb Z_p\hoge[G]
\simeq
\mathbb Z_p\hoge[T_1,\ldots,T_d],
\qquad
\gamma_i\longmapsto 1+T_i,
\]
or equivalently $\gamma_i-1\mapsto T_i$.

\begin{defin}
A finitely generated $\Lambda_G$-module $N$ is called a pseudo-null module
if there exist coprime elements $\phi,\psi\in\Lambda_G$ such that
$\phi,\psi\in\operatorname{Ann}_{\Lambda_G}(N)$.
\end{defin}

For a finitely generated $\Lambda_G$-module $M$, we use the following
structure theorem.
\begin{lemma}[\cite{Wa} Section 13.2]\label{fg-lambda}
Let $M$ be a finitely generated $\Lambda_G$-module. Then there exists an
exact sequence
\[
0\to Z\to M\to
\Lambda_G^a\oplus
\bigoplus_{i=1}^{\lambda(M)}\Lambda_G/\mathfrak p_i^{e_i}
\to Z'\to 0,
\]
where $Z$ and $Z'$ are pseudo-null $\Lambda_G$-modules and the
$\mathfrak p_i$ are height-one prime ideals of $\Lambda_G$.
\end{lemma}

The integer $a$ in Lemma \ref{fg-lambda} is called the $\Lambda_G$-rank
of $M$, and is denoted by $\operatorname{rk}_{\Lambda_G}(M)$. The
characteristic ideal of the torsion part of $M$ is defined by
$\operatorname{char}_{\Lambda_G}(M):=\prod_{i=1}^{\lambda(M)}
\mathfrak p_i^{e_i}$.

If there exists a morphism between finitely generated $\Lambda_G$-modules
$M\to M'$ whose kernel and cokernel are pseudo-null, we say that $M$ and
$M'$ are pseudo-isomorphic.

Greenberg \cite[Theorem 1]{Gr} proved the following result.
\begin{prop}\label{General-Greenberg}
Let $k$ be a number field, and let $S$ be a finite set of places of $k$ above
$p$. Let $k_\infty/k$ be a $\mathbb Z_p^d$-extension unramified outside $S$,
and let $\mathfrak{X}_{k,S}$ be the Galois group of the maximal abelian pro-$p$ extension
of $k_\infty$ unramified outside $S$.
\begin{enumerate}
\item Suppose $d=1$ and that the natural localization map
$\varprojlim_n\mathcal O_{k_n}^{\times}\otimes\mathbb Z_p\to
\varprojlim_n\prod_{v\in S}\mathcal O_{k_n,v}^{\times}\widehat{\otimes}\mathbb Z_p$
is injective, where $k_n$ is the $n$-th layer of $k_\infty/k$. Then $\mathfrak{X}_{k,S}$
contains no non-trivial finite submodule.

\item For general $d$, assume that the corresponding localization map
$\mathcal O_{k}^{\times}\otimes\mathbb Z_p\to
\prod_{v\in S}\mathcal O_{k,v}^{\times}\widehat{\otimes}\mathbb Z_p$
is injective. Then $\mathfrak{X}_{k,S}$ contains no non-trivial pseudo-null
$\Lambda$-submodule.
\end{enumerate}
\end{prop}

\begin{remark}
\begin{enumerate}\label{remark:coh-leopold}
\item
Maire proves in \cite[Proposition 3.1]{Maire2005} that the vanishing of 
$H^2(G_{k}^{S},\mathbb Q_p/\mathbb Z_p)$ holds when the following localization map is injective;
\[\mathcal O_{k}^{\times}\widehat{\otimes}\mathbb Z_p
\longrightarrow
\prod_{v\in S}\mathcal O_{k,v}^{\times}\widehat{\otimes}\mathbb Z_p.
\]
This condition is known in several important cases. Brumer proved such
injectivity when $k$ is an abelian extension of an imaginary quadratic field
$k_0$ and $S$ is the set of primes above a fixed prime $\mathfrak p\mid p$
of $k_0$. 

\item
The theorem also holds for general finite sets of finite places $S$, with the
same localization condition. If real places are included, however, the
statement has to be modified.
\end{enumerate}
\end{remark}

The following theorem is standard, but we include a proof because we could not
find a reference in the precise form needed below.
\begin{theorem}\label{thm:projective-resolution}
Let $k$ be a number field and let $k_\infty/k$ be a 
$\Z_p^r$-extension unramified outside $S$ where $r=1$ or $2$. Let
$\Lambda:=\mathbb Z_p\hoge[\operatorname{Gal}(k_\infty/k)]$.
Assume that the corresponding localization map
$\mathcal O_{k}^{\times}\otimes\mathbb Z_p\to
\prod_{v\in S}\mathcal O_{k,v}^{\times}\widehat{\otimes}\mathbb Z_p$
is injective. Assume moreover that
$\mathfrak{X}_{k,S}\neq 0$ and that one of the following conditions holds:
\begin{enumerate}
\item There exists a $\mathbb Z_p$-extension inside $k_\infty/k$, unramified
outside $S$, whose corresponding Iwasawa module has trivial $\mu$-invariant.
\item The set $S$ contains all primes of $k$ above $p$.
\end{enumerate}
Then $\mathfrak{X}_{k,S}$ contains no non-trivial finite submodule and admits a projective
resolution
\[
0\to \Lambda^{d-1-s}\to \Lambda^{d-1}\to \mathfrak{X}_{k,S}\to 0,
\]
where $s=\mathrm{rank}_\Lambda \mathfrak{X}_{k,S}$, and \(
d:=\dim_{\mathbb F_p}H^1(G_k^{S},\mathbb F_p),
\).
\end{theorem}
\begin{proof}
We first note that $H^2(G_k^{S}, \Q_p/\Z_p)=0$ by Remark \ref{remark:coh-leopold} (1). We also note that the vanishing of $H^2$ descends to any 
$\mathbb Z_p$-extensions considered below. Let $K/k$ be a $\mathbb Z_p$-
extension of $k$ unramified outside $S$,
and put $\Gamma:=\operatorname{Gal}(K/k)$. We have an exact sequence
\[
1\longrightarrow G_K^S\longrightarrow G_k^S\longrightarrow \Gamma
\longrightarrow 1.
\]
The Hochschild--Serre spectral sequence
\[
H^i(\Gamma,H^j(G_K^S,\mathbb Q_p/\mathbb Z_p))
\Rightarrow H^{i+j}(G_k^S,\mathbb Q_p/\mathbb Z_p)
\]
shows that $H^2(G_K^S,\mathbb Q_p/\mathbb Z_p)^\Gamma$ occurs as a graded
quotient of $H^2(G_k^S,\mathbb Q_p/\mathbb Z_p)$, because
$\operatorname{cd}_p(\Gamma)=1$. Hence, if
$H^2(G_k^S,\mathbb Q_p/\mathbb Z_p)=0$, then
$H^2(G_K^S,\mathbb Q_p/\mathbb Z_p)^\Gamma=0$. Since
$H^2(G_K^S,\mathbb Q_p/\mathbb Z_p)$ is a discrete $p$-primary
$\Gamma$-module, the vanishing of its $\Gamma$-invariants implies that the
module itself is zero. Therefore
\[
H^2(G_K^S,\mathbb Q_p/\mathbb Z_p)=0.
\]
If $k_\infty/k$ is a $\mathbb Z_p^2$-extension, we choose an intermediate
field $K$ such that $K/k$ and $k_\infty/K$ are both $\mathbb Z_p$-extensions,
and apply the preceding argument twice. Thus
$H^2(G_{k_\infty}^S,\mathbb Q_p/\mathbb Z_p)=0$.
Next, we will recall the result of Nguyen-Quan-Do \cite{NguyenQuangDo1980RelationsPGroupe}.
Let $G:=G_{k}^{S}$ be the Galois group of the maximal pro-$p$ extension for $k$ unramified outside $S$, and put $d:=d(G)=\dim_{\F_p}H^1(G,\F_p)$, the minimal number of topological generators of $G$. Choose a free pro-$p$ group $F=\kakko<x_1,\ldots,x_d>$ and a minimal continuous presentation $F\twoheadrightarrow G$.
Let $\Gamma:=\Gal(k_\infty/k)$ and let 
$\pi:F\rightarrow G\rightarrow \Gamma$ be the induced map. Put $R:=\ker(F\rightarrow \Gamma)$ and $N:=\ker(F\rightarrow G)$.
We also put $\Lambda_F:=\Z_p\hoge[F]$ and $\Lambda_\Gamma:=\Z_p\hoge[\Gamma]$, and let $I_F$ be the augmentation ideal of $\Lambda_F$. 

Let $g_i:=\pi(x_i)\in\Gamma$. We use the same symbol for an element of
$\Lambda_F$ and for its image in $\Lambda_\Gamma$ when no confusion can
arise. Similarly, the Fox derivative $\partial/\partial x_i$ is first defined
on $\Lambda_F$ and then composed with the natural map
$\Lambda_F\to\Lambda_\Gamma$.

The Fox derivatives are characterized by
$\partial x_i/\partial x_j=\delta_{ij}$ and
\[
\frac{\partial uv}{\partial x_i}
=
\frac{\partial u}{\partial x_i}
+u\frac{\partial v}{\partial x_i}
\]
for $u,v\in F$. The fundamental identity for Fox derivatives says that, for
any $w\in F$,
\[
w-1=\sum_{i=1}^d\frac{\partial w}{\partial x_i}(x_i-1)
\quad\text{in } I_F.
\]
Applying $\Lambda_F\to\Lambda_\Gamma$, we get
\[
\pi(w)-1
=
\sum_{i=1}^d
\pi\!\left(\frac{\partial w}{\partial x_i}\right)(g_i-1).
\]
In particular, if $n\in R$, then $\pi(n)=1$, and hence
$\sum_{i=1}^d(\partial n/\partial x_i)(g_i-1)=0$.

Thus the map
\[
R^{\mathrm{ab}}\to\Lambda_\Gamma^d,\qquad
n\mapsto \left(\frac{\partial n}{\partial x_i}\right)_{i=1}^d
\]
has image contained in the kernel of
$\Lambda_\Gamma^d\to I_\Gamma$, $e_i\mapsto g_i-1$. Lyndon \cite{Lyndon50} proved that this image is exactly the kernel. Therefore we
obtain a short exact sequence 
\begin{equation}\label{eq:Lydon-short-exact}
0\to R^{\mathrm{ab}}\to\Lambda_\Gamma^d\to I_\Gamma\to 0.
\end{equation}
Here, $R^{\mathrm{ab}}$ is a $\Gamma$-module by lifting $\gamma\in \Gamma$ to $\widetilde{\gamma}\in F$ and taking conjugation.
Then, we have $\partial \widetilde{\gamma}n\widetilde{\gamma}^{-1}/\partial x_i=\partial \widetilde{\gamma}/\partial x_i+\pi(\widetilde{\gamma})\partial n/\partial x_i-\pi(\widetilde{\gamma}n\widetilde{\gamma}^{-1})\partial \widetilde{\gamma}/\partial x_i=\gamma \partial n/\partial x_i$ and hence, the short exact sequence (\ref{eq:Lydon-short-exact}) is a $\Lambda$-module exact sequence. Applying the same argument with \(N:=\mathrm{ker}(F\to G)\), \(\Lambda_G:=\mathbb Z_p\hoge[G]\), and the augmentation ideal \(I_G\subset \Lambda_G\), we obtain a short exact sequence of \(\Lambda_G\)-modules
\begin{equation}\label{eq:Lydon-2}
0\rightarrow N^{\mathrm{ab}} \rightarrow \Lambda_G^d\rightarrow I_G\rightarrow 0
\end{equation}
Let
\(
H:=G_{k_\infty}^{S}=\ker(G\to\Gamma).
\)
Taking \(H\)-homology in the exact sequence \((\ref{eq:Lydon-2})\), we obtain
\begin{equation}\label{eq:Nguen}
\ldots \rightarrow H_1(H, I_G)\rightarrow N^{\mathrm{ab}}_{H}\rightarrow \Lambda^d\rightarrow (I_G)_H\rightarrow 0.
\end{equation}
Since $\Z_p\hoge[G]$ is an induced $H$-module, Shapiro's lemma gives $H_i(H, \Z_p\hoge[G])=0$ for $i>0$. From the $H$-module exact sequence $0\rightarrow I_G\rightarrow \Z_p\hoge[G]\rightarrow \Z_p\rightarrow 0$, we have 
$H_2(H,\Z_p)\simeq H_1(H,I_G)$, and $0\rightarrow H_1(H,\Z_p) \rightarrow (I_G)_H\rightarrow \Lambda\rightarrow \Z_p\rightarrow  0$. The first isomorphism implies
\begin{equation}\label{eq:aa}
H_1(H, I_G)\simeq H_2(G_{k_\infty}^S,\Z_p)\simeq H^2(G_{k_\infty}^S,\Q_p/\Z_p)^{\vee}=0,
\end{equation}
and the second exact sequence implies 
\begin{equation}\label{eq:bb}
0\rightarrow \mathfrak{X}_{k,S}\rightarrow (I_G)_H\rightarrow I_\Gamma \rightarrow 0.
\end{equation}
Combined with (\ref{eq:Nguen}), (\ref{eq:aa}), we obtain the following commutative diagram;
\[
\begin{CD}
0 @>>>N^{\mathrm{ab}}_{H}@>>> \Lambda^d@>>> (I_G)_H@>>> 0\\
@. @VVV @VVV @VVV @.\\
0 @>>> 0 @>>>I_\Gamma @>>>I_\Gamma @>>>0
\end{CD}\]
The middle vertical morphism is from (\ref{eq:Lydon-short-exact}), and the right down side morphism is from (\ref{eq:bb}).
Then, by the snake lemma, we have 
\[
0 \rightarrow  N^{\mathrm{ab}}_H \rightarrow R^{\mathrm{ab}} \rightarrow \mathfrak{X}_{k,S} \rightarrow 0.\]
By \cite[Proposition~3.4]{Maire2005} and \cite[Theorem 5.1]{Uwe}, under either assumption (1) or (2),
the pro-$p$ group $G_k^S(p)$ has $p$-cohomological dimension at most $2$.
Therefore, Jannsen's theorem \cite[Lemma $4.3$ (b)]{Uwe} applies, and
\[
(N^{\mathrm{ab}})_H=N/[N,R]
\]
is a free $\Lambda_\Gamma$-module. Next, we will show 
$R^{\mathrm{ab}}\simeq\Lambda_\Gamma^{d-1}$ by using
\eqref{eq:Lydon-short-exact}.
Since $\Gamma\simeq \mathbb Z_p^r$ with $r=1$ or $2$, we may choose the
minimal generators $x_1,\ldots,x_d$ of $F$ so that $x_1,\ldots,x_r$ map to
topological generators $\gamma_1,\ldots,\gamma_r$ of $\Gamma$ and
$x_{r+1},\ldots,x_d$ map to $1$. Hence the map
$\Lambda_\Gamma^d\to I_\Gamma$ is given by
\[
e_i\mapsto \gamma_i-1\quad (1\leq i\leq r),\qquad
e_i\mapsto 0\quad (i>r).
\]
For $r=1$ its kernel is freely generated by $e_2,\ldots,e_d$. For $r=2$ its
kernel is freely generated by
\[
(\gamma_2-1)e_1-(\gamma_1-1)e_2,\ e_3,\ldots,e_d.
\]
Thus in both cases
\[
R^{\mathrm{ab}}\simeq \ker(\Lambda_\Gamma^d\to I_\Gamma)
\simeq \Lambda_\Gamma^{d-1}.
\]
Hence, we obtain a projective resolution with some $s\ge 0$,
\[
0\to \Lambda_\Gamma^{d-1-s}\to \Lambda_\Gamma^{d-1}\to \mathfrak{X}_{k,S}\to 0.
\] 

This proves the theorem.
\end{proof}
From this proof, we also obtain the following corollary.

\begin{corollary}\label{cor:chosen-generators-resolution}
Let $x_1,\ldots,x_{d-1}$ be a chosen generating system of $\mathfrak{X}_{k,S}$ as a
$\Lambda$-module. Then the surjection
$\pi:\Lambda^{d-1}\to \mathfrak{X}_{k,S}$ in Theorem \ref{thm:projective-resolution} may
be chosen so that, for the standard $\Lambda$-basis
$e_1,\ldots,e_{d-1}$ of $\Lambda^{d-1}$, one has
$\pi(e_i)=x_i$ for $1\leq i\leq d-1$.
\end{corollary}

\section{Topological Generators}

B\o{}lling \cite{Bo} derived a topological generating system for the Galois group $G_K^{(S)}(p)$ over $K$ for an arbitrary algebraic number field $K$ and a finite set of primes $S$, where $G_K^{(S)}(p)$ is the Galois group of the maximal pro-$p$ extension of $K$ unramified outside $S$.  

In this section, based on B\o{}lling's method, we derive the topological generating system of $G^{(\pi)}_K$ for an imaginary quadratic field $K$. Put $G^\pi:=G_K^\pi$ for simplicity. We first derive a fundamental short exact sequence from class field theory. Let $C_K = \A_K^\times / K^\times$ denote the idele class group, and set  
\[
U_K(\pi) := \prod_{v \nmid \infty \pi} \O_{K_v}^\times \times \C \subset \A_K^\times.
\]  
Let $\overline{U_K(\pi)}$ be the closure of the image of $U_K(\pi)$ in $\A_K^\times \to C_K = \A_K^\times / K^\times$. Similarly, define $U_K := \prod_{v \nmid \infty} \O_{K_v}^\times \times \C$ and let $\overline{U_K} \subset C_K$ be its closure.

Class field theory provides the following exact sequence:
\[
1 \to \overline{U_K} \to C_K \to \mathrm{Cl}_K \to 1.
\]  
Since $\overline{U_K(\pi)} \subset \overline{U_K}$, it follows that
\[
1 \to \overline{U_K}/\overline{U_K(\pi)} \to C_K/\overline{U_K(\pi)} \to \mathrm{Cl}_K \to 1.
\]  
Moreover, class field theory implies that if $K_\pi^{\mathrm{ab}}$ is the maximal abelian extension of $K$ unramified outside $\pi$, then
\[
C_K / \overline{U_K(\pi)} \simeq \mathrm{Gal}(K_\pi^{\mathrm{ab}}/K).
\]  
Taking the maximal pro-$p$ quotient, equivalently applying completed tensor product $\widehat{\otimes} \Z_p$ to this short exact sequence yields the following:
  
\begin{lemma}\label{fundamental}
Let $A_K := \mathrm{Cl}_K \otimes \Z_p$. Then, there is a short exact sequence:
\[
1 \to (\overline{U_K} / \overline{U_K(\pi)}) \widehat{\otimes} \Z_p \to (G^{\pi})^{\mathrm{ab}} \to A_K \to 1,
\]
where $G^{\pi}$ is the Galois group of the maximal pro-$p$ extension of $K$ unramified outside $\pi$.
\end{lemma}

Assuming that $\overline{\O_{K,\pi}^\times \cap K^\times U_K(\pi)} \subset \O_{K,\pi}$ denotes the topological closure, it follows that
\[
\overline{U_K} / \overline{U_K(\pi)} \widehat{\otimes} \Z_p \simeq \O_{K,\pi}^\times / \overline{\O_{K,\pi}^\times \cap K^\times U_K(\pi)} \widehat{\otimes} \Z_p.
\]  
According to Washington \cite{Wa}, $13.1, p.266$, $\overline{\O_{K,\pi}^\times \cap K^\times U_K(\pi)}$ can be described as the topological closure of the image of
\[
\iota: (1 + \pi \O_K) \cap \O_K^\times \to 1 + \pi \O_{K,\pi} \subset \O_{K,\pi}^\times.
\]  
Since $\O_K^\times$ is a finite group, the topological closure of this image is also finite. Thus, the following lemma follows from Lemma \ref{fundamental}:

\begin{prop}\label{Fundamental}
Let $K$ be an imaginary quadratic field. Then, there is a short exact sequence:
\[
1 \to \Z_p \to (G^{\pi})^{\mathrm{ab}} \to A_K \to 1.
\]
\end{prop}

\begin{proof}
Since $(\overline{U_K} / \overline{U_K(\pi)}) \widehat{\otimes} \Z_p \simeq (\O_{K,\pi}^\times / \overline{\O_{K,\pi}^\times \cap K^\times U_K(\pi)}) \widehat{\otimes} \Z_p$ and we have already shown that 

$\overline{\O_{K,\pi}^\times \cap K^\times U_K(\pi)}$ is a finite group, it follows that $\O_{K,\pi}^\times \otimes \Z_p \simeq \Z_p$, and the finite part becomes trivial. Therefore,
\[
(\overline{U_K} / \overline{U_K(\pi)}) \widehat{\otimes} \Z_p \simeq \Z_p.
\]
\end{proof}

Using this proposition, we provide a topological generating system for $G^{\pi}$. We shall use the following Burnside basis theorem for pro-$p$ groups (see, for example, Koch \cite[Theorem 4.2]{Ko}).

\begin{prop}\label{Burnside}
Let $G$ be a pro-$p$ group. Then the minimal number of topological generators of $G$ is
\[
\dim_{\F_p}G/\overline{[G,G]}G^p=\dim_{\F_p}H^1(G,\F_p).
\]
Equivalently, any lift to $G$ of an $\F_p$-basis of $G/\overline{[G,G]}G^p$ is a minimal system of topological generators of $G$.
\end{prop}

We now prove the following proposition. Using the map $(G^{\pi})^{\mathrm{ab}} \to A_K$ in Lemma \ref{fundamental}, define
\[
B_K := \text{image of } ((G^{\pi})^{\mathrm{ab}}_{\mathrm{tors}}\rightarrow A_K).
\]

\begin{prop}\label{Iwasawa-module-top}
Let $K_\infty$ be the maximal $\Z_p$-extension of an imaginary quadratic field $K$ unramified outside $\pi$, and let $\Gamma := \Gal(K_\infty/K)$. For the Hilbert class field $K(1)$ of $K$, let $p^n := |[K_\infty \cap K(1) : K]|$. Then,
\(
|B_K| = |A_K|/p^n.
\)
Let $a_1, \ldots, a_r$ be the generators of $B_K$. Then the following holds:
\begin{enumerate}
    \item $B_K \simeq (G^{\pi})^{\mathrm{ab}}_{\mathrm{tors}}$.
    \item We identify $(G^{\pi})^{\mathrm{ab}}$ with $\Gamma \times B_K$. The lift to $G^{\pi}$ of the topological generators $\gamma_0, a_1, \ldots, a_r$ of $\Gamma \times B_K$ gives the minimal topological generating system $\gamma, g_1, \ldots, g_r$ for $G^{\pi}$.
    \item $\mathfrak{X}_{K,\pi} = \Gal(M_\pi(K)/K_\infty)$ satisfies
    \[
    \mathfrak{X}_{K,\pi} = \left(\langle g_1, \ldots, g_r \rangle \overline{[G^{\pi}, G^{\pi}]}\right)^{\mathrm{ab}}.
    \]
    Furthermore, when $r = 1$, $\mathfrak{X}_{K,\pi}$ is cyclic as a $\Lambda_\Gamma$-module, that is,
    \[
    \mathfrak{X}_{K,\pi} = \Lambda_\Gamma \cdot g_1.
    \]
\end{enumerate}
\end{prop}

\begin{proof}
Applying $\operatorname{Tor}_1^{\Z_p}(-,\Q_p/\Z_p)$ to Proposition~\ref{Fundamental}, we obtain the following short exact sequence:
\[
1 \to (G^{\pi})^{\mathrm{ab}}_{\mathrm{tors}} \to (A_{K})_{\mathrm{tors}} = A_{K}.
\]
By the definition of $B_K$, the image of the map to $A_K$ in this sequence is $B_K$. Thus, we obtain $(G^{\pi})^{\mathrm{ab}}_{\mathrm{tors}} \simeq B_K$. 

Next, we prove the second claim. By Proposition~\ref{Fundamental} and the theory of elementary divisors, we have
\[
(G^{\pi})^{\mathrm{ab}} \simeq \Z_p \times ((G^{\pi})^{\mathrm{ab}})_{\mathrm{tors}}.
\]
By the first claim,
\begin{equation}\label{torsion+free}
((G^{\pi})^{\mathrm{ab}})_{\mathrm{tors}} \simeq B_K.
\end{equation}
Using Burnside's basis theorem (Proposition~\ref{Burnside}), the minimal topological generating system for $G^{\pi}$ is given by the lift to $G^{\pi}$ of the minimal topological generating system for $(G^{\pi})^{\mathrm{ab}}$. By equation (\ref{torsion+free}), $G^{\pi}$ is generated by the lift of the topological generator of $\Gamma \simeq \Z_p$, denoted by $\gamma$, and the lift of the minimal generating set for $B_K$, denoted by $g_1, \ldots, g_r$. 

Finally, we prove the third claim. For the field extension $K_\pi(p)/M_\pi(K)/K_\infty$, we have the following natural short exact sequence:
\[
1 \to [G^\pi, G^\pi] \to \Gal(K^{\pi}(p)/K_\infty) \to (G^{\pi})^{\mathrm{ab}}_{\mathrm{tors}} \to 1.
\]
Thus, if $\tau: (G^{\pi})^{\mathrm{ab}}_{\mathrm{tors}} \to \Gal(K^\pi(p)/K_\infty)$ is a set-theoretic section, then 
\[
\Gal(K^\pi(p)/K_\infty) = \kakko<\tau((G^\pi)^{\mathrm{ab}}_{\mathrm{tors}})>[G^\pi, G^\pi],
\]
and we obtain $\mathfrak{X}_{K,\pi} = (\langle g_1, \ldots, g_r \rangle [G^\pi, G^\pi])^{\mathrm{ab}}$. 

In particular, when $r = 1$, $G^\pi = \langle \gamma, g_1 \rangle$. Let $\overline{\kakko<g_1>}=\Lambda_\Gamma\cdot g_1$ be the normal closed subgroup of $G^\pi$ containing $g$. Then, each element $x,y\in G^\pi$ can be written as $x=x_1y_1$, $y=x_2y_2$ where $y_1,y_2\in \overline{\kakko<g_1>}$, and $x_1,x_2\in \kakko<\gamma>$.  Hence, we have
\[
[x,y] = xyx^{-1}y^{-1} = y_1^{x_1}(y_2y_1^{-1})^{x_1x_2}(y_2^{-1})^{x_2} \in \Lambda_\Gamma \cdot g_1.
\]
Here, for $y\in \overline{\kakko<g_1>}$, $x\in \kakko<\gamma>$, we denoted $y^x:=xyx^{-1}$.
Thus, in $\mathfrak{X}_{K,\pi}$, we obtain $(\langle g_1 \rangle [G^\pi, G^{\pi}])^{\mathrm{ab}} \subset \Lambda_\Gamma \cdot g_1$. Since $g_1 \in \mathfrak{X}_{K,\pi}$ and $\Lambda_\Gamma \cdot g_1 \subset \mathfrak{X}_{K,\pi}$, we have $\mathfrak{X}_{K,\pi} = \Lambda_\Gamma \cdot g_1$.
\end{proof}
When $K_\infty$ is the maximal $\Z_p^2$-extension, we obtain the following proposition.
\begin{prop}
Let $K_\infty/K$ be the maximal $\Z_p^2$-extension, and let $K^{\mathrm{anti}}$ be its maximal anticyclotomic $\Z_p$-extension. Assume that $\mathrm{Cl}_K\otimes \Z_p$ is cyclic and $K^{\mathrm{anti}}$ does not contain the maximal unramified $p$-extension of $K$.
Assume moreover that either $p$ splits in $K$, or $p\neq 3$, or $\mathrm{disc}(K)\not\equiv -3 \mod 9$. Let $G^{p}=\Gal(K^{(p)}/K)$ be the Galois group of the maximal pro-$p$ extension of $K$ unramified outside $p$. Then \(G^{p}\) has a minimal system of topological generators
\(
G^{p}=\langle \gamma_1,\gamma_2,g\rangle,
\)
where \(\gamma_1,\gamma_2\) are lifts of topological generators of
\(\operatorname{Gal}(K_\infty/K)\), and
\[
\mathfrak X_{K,p}
=
\Lambda\cdot \pi(g)
+
\Lambda\cdot \pi([\gamma_1,\gamma_2]).
\]
Here, $\mathfrak{X}_{K,p}$ is the Galois group of the maximal abelian pro-$p$ extension of $K_\infty$ unramified outside $p$ and $\pi:\Gal(K^{(p)}/K_\infty)\rightarrow \mathfrak{X}_{K,p}$ is the canonical restriction map.
\end{prop}
\begin{remark}
By class field theory, the maximal unramified extension $H/K$ is anticyclotomic. In other words, for all $g\in \Gal(H/K)$, we have $cgc^{-1}=g^{-1}$ with the complex conjugate $c\in \Gal(K/\Q)$. 
\end{remark}
\begin{proof}
For a number field $L$, we denote by $M(L)$ the maximal abelian $p$-extension of $L$ unramified outside $p$.
First, we discuss topological generators. We study only the case when $p$ splits in $K/\Q$ because we can prove the remaining case with the same arguments. $G^p$ is topologically generated by the lift of the generators of $(G^p)^{\mathrm{ab}}=\Gal(M(K)/K)$, where $M(K)$ is the maximal abelian $p$-extension of $K$ unramified outside $p$. By class field theory, we obtain the exact sequence
$1\rightarrow \Z_p\times \Z_p\rightarrow \Gal(M(K)/K)\rightarrow \mathrm{Cl}_K\otimes \Z_p\rightarrow 1$. Hence, we have $G^p=\kakko<g,\gamma_1,\gamma_2>$, where $\gamma_1,\gamma_2$ generate $\Gal(K_\infty/K)$ and $g\in \Gal(K^{(p)}/K_\infty)$.

Let $\kakko<g,\gamma_1,\gamma_2>$ be topological generators of $G^p$. Let $\widetilde{\Gamma}$ be the subgroup of $G^p$ topologically generated by $\gamma_1,\gamma_2$ and $\overline{\kakko<g>}$ be the minimal normal closed subgroup of $G^p$ containing $g$.  

Step $1$. 
We will show that $\Gal(K^{(p)}/K_\infty)$ is topologically generated by 
$[G^{p},G^p]$ and $g$. 

Consider the exact sequence $1\rightarrow \Gal(M(K_\infty)/M(K))=[G^p, G^p]^{\mathrm{ab}}\rightarrow \Gal(M(K_\infty)/K_\infty)=\mathfrak{X}_{K,p}\rightarrow \Gal(M(K)/K_\infty)\rightarrow 1$, and obtain the claim.

Step $2$. we prove that $\Gal(K^{(p)}/K_\infty)$ is topologically generated by $[\widetilde{\Gamma},\widetilde{\Gamma}]$, and $\overline{\kakko<g>}$. 

$\overline{\kakko<g>}$ is a normal subgroup of $G^{p}=\Gal(K^{(p)}/K)$ and hence, we know that $G^p=\widetilde{\Gamma}\cdot \overline{\kakko<g>}$. Now, for $x_1,x_2\in \widetilde{\Gamma}$, and $y_1,y_2\in \overline{\kakko<g>}$, we have $[x_1y_1,x_2y_2]= xyx^{-1}y^{-1} = y_1^{x_1}(y_2y_1^{-1})^{x_1x_2}[x_1,x_2](y_2^{-1})^{x_2}$ holds. Here we denote $y^{x}:=xyx^{-1}$ for $x\in \widetilde{\Gamma}$, $y\in G^p$. 
Therefore, we conclude that $[G^p,G^p]$ is contained in the subgroup of $\Gal(K^{(p)}/K_\infty)$ generated by $[\widetilde{\Gamma},\widetilde{\Gamma}]$ and $\overline{\kakko<g>}$. Combined with Step $1$, we obtain that $\Gal(K^{(p)}/K_\infty)$ is generated by $[\widetilde{\Gamma},\widetilde{\Gamma}]$ and $\overline{\kakko<g>}$.

Step 3. we prove that $\mathfrak{X}_{K,p}=\Lambda\cdot \pi([\gamma_1,\gamma_2])+\Lambda\cdot \pi(g)$ holds, where $\pi:\Gal(K^p/K_\infty)\twoheadrightarrow \mathfrak{X}_{K,p}$ is the canonical restriction map. 

In general, for a free pro-$p$ group generated by two elements $B=\kakko<x,y>$, $B$ acts on $[B,B]^{\mathrm{ab}}$ by conjugation and as a $\Z_p\hoge[B^{\mathrm{ab}}]$-module, $[B,B]^\mathrm{ab}$ is a cyclic $\Z_p\hoge[B^{\mathrm{ab}}]$-module (see, Ihara\cite{Iha} Theorem $2$ (i)).
Therefore, we know that $[\widetilde{\Gamma},\widetilde{\Gamma}]^{\mathrm{ab}}$ is a cyclic $\Z_p\hoge[\Gamma]$-module, and hence $[\widetilde{\Gamma},\widetilde{\Gamma}]^{\mathrm{ab}}=\Lambda\cdot \pi([\gamma_1,\gamma_2])$. By definition, $\pi(\overline{\kakko<g>})=\Lambda\cdot g$ and hence, we obtain our result.
\end{proof}
\section{Structure of Iwasawa Modules}
In this section, we recall the Iwasawa main conjecture of an imaginary quadratic field proved by K. Rubin \cite{Rubin91}, \cite{Rubin1994MoreMainConjectures}.
Let $K$ be an imaginary quadratic field, and let $K_\infty$ be the maximal $\mathbb{Z}_p$-extension of $K$ unramified outside $\pi$ or the maximal $\Z_p^2$-extension of $K$. In the following, we assume that $p$ is an odd prime and splits in $K/\Q$.

First, we introduce a $p$-adic $L$-function. We fix embeddings $\iota_\infty: \overline{\mathbb{Q}} \hookrightarrow \mathbb{C}$ and $\iota_p: \overline{\mathbb{Q}} \hookrightarrow \mathbb{C}_p$ such that $\iota_p^{-1}(m_{\mathbb{C}_p}) \cap \mathcal{O}_K = \pi$, where $m_{\mathbb{C}_p}$ is the maximal ideal of $\mathcal{O}_{\mathbb{C}_p}$. 

We also fix an embedding $\sigma:K\hookrightarrow \overline{\mathbb{Q}}$. Let $\chi: I_K^{\mathfrak{f}} \to \overline{\mathbb{Q}}^\times$ be an algebraic Hecke character of conductor $\mathfrak{f}$ with infinite type $(-k, 0)$ for some $k \geq 1$, where $I_K^{\mathfrak{f}}$ denotes the group of fractional ideals of $K$ coprime to $\mathfrak{f}$. In other words, $\chi$ is a character such that for any $a \in \mathcal{O}_K$ satisfying $a \equiv 1 \mod \mathfrak{f}$, we have $\chi(a) = \sigma(a)^{-k}$. We always use the notation $\chi_{\mathrm{gal}}:G_K^{(\pi\mathfrak{f})}\rightarrow \O_{\C_p}^\times$ to express the Galois representation attached to $\chi$.
For such a character $\chi$, consider the $L$-function
\[
L_\mathfrak{f}(s, \chi) := \sum_{\mathfrak{a} \text{ integral, coprime to } \mathfrak{f}} \frac{\iota_\infty \circ \chi(\mathfrak{a})}{N(\mathfrak{a})^s},
\]
which converges absolutely for $\mathrm{Re}(s) > 1 + k/2$. This function can be analytically continued to the entire complex plane if $\chi$ is not the power of the Norm character. We define $L_\mathfrak{f}(\chi) := L_\mathfrak{f}(0, \chi)\in \C$.
 
Let \(E/K(1)\) be an elliptic curve defined over the Hilbert class field such that
\(
\iota:\operatorname{End}_{K(1)}(E)\simeq \mathcal O_K,
\)
and suppose that \(E\) has  good reduction at the prime $v=\iota_p^{-1}(m_{\C_p})\cap \O_{K(1)}$ above $\pi$. We choose the isomorphism $\iota$ such that $\omega\circ \iota(a)=a\omega$ for all $a\in \O_K$, and $\omega\in \Gamma(E_{/K(1)},\Omega^1_{E/K(1)})$. 
\begin{defin}\label{def:period}
Let $\mathcal{E}_{\O_{(v)}}$ be the Neron model of $E_{/K(1)}$ where $\O_{(v)}:=\O_{K(1),(v)}$. Choose $\omega \in \Gamma(\mathcal{E}, \Omega^1_{\mathcal{E}/\O_{(v)}})$ to be a free $\O_v$-basis of $\Gamma(\mathcal{E}, \Omega^1_{\mathcal{E}/\O_{(v)}})$. Let $\widehat{\mathcal{E}}$ be the formalization of $\mathcal{E}$ along the unit section and denote $\widehat{\omega}\in \Gamma(\widehat{\mathcal{E}},\widehat{\Omega}^1_{\widehat{\mathcal{E}}/\O_{v}})$ be the differential induced by $\omega$. We use the notation $W:=W(\bar\F_p)$. Now, fix two isomorphisms 
$$\theta:\C/\O_K\simeq E(\C)\ ,\ \theta_p:\widehat{\mathbb{G}}_{m/W}\simeq \widehat{\mathcal{E}}_{/W}.$$
The existence of the second isomorphism is proved as follows: Let $\overline{\mathcal{E}}$ be the reduction of $\mathcal{E}$. As a finite group scheme, we have $\overline{\mathcal{E}}[p^n] \simeq \Z/p^n\Z$ over $\bar\F_p$, and hence, we have $\mathcal{E}[p^\infty]^{\mathrm{et}}\simeq \Q_p/\Z_p$ over $W$ by the existence of Serre-Tate's canonical lift. Therefore, we have an isomorphism of $p$-divisible group $\mathcal{E}[p^\infty]^\circ\simeq \mu_{p^\infty}$ over $W$ due to Cartier duality.

Choosing such isomorphisms is equivalent to choosing an \(\mathcal O_K\)-basis of \(H_1(E(\mathbb C),\mathbb Z)\) and a \(\mathbb Z_p\)-basis of
\(
\varprojlim_n \mathcal E[p^n]^\circ(\mathcal O_{\mathbb C_p}).
\)
Define $\Omega\in \C^\times$, $\Omega_p\in W^\times$ as follows;
$$\theta^{*}\omega=\Omega dz\ ,\  \theta_p^{*}\widehat{\omega}=\Omega_p\widehat{\frac{dt}{t}}.$$
\end{defin}
We remark that this pair $(\Omega, \Omega_p)\in \C^\times \times W^\times$ depends on the choice of $(E, \omega,\theta, \theta_p)$. 

\begin{theorem}[( \cite{ds}, Theorem 4.12)]\label{p-adic-interpolate}
Let $K$ be an imaginary quadratic field, and let $\mathfrak{f}$ be an integral ideal coprime to $\pi$. Fix the pair $(E,\omega, \theta,\theta_p)$ as in Definition \ref{def:period} and define $(\Omega,\Omega_p)$. Let $\mathcal{G}_\mathfrak{f}:=\Gal(K(\mathfrak{f}\pi^\infty)/K)$.
\begin{enumerate}
  \item[$(1)$.] There exists a unique $p$-adic measure $\mu_\mathfrak{f} \in \Lambda_{\mathcal{G}_\mathfrak{f}} \otimes_{\mathbb{Z}_p} W$ which depends only on the choice of $(E,\theta,\theta_p)$ and satisfies
  \[
  \Omega_p^{-k} \int_{\mathcal{G}_\mathfrak{f}} \chi_{\mathrm{gal}}^{-1}(\sigma) d\mu_\mathfrak{f}(\sigma) = \Omega^{-k} \cdot G(\chi) \cdot \left(1 - \frac{\chi(\pi)}{p}\right) \cdot L_\mathfrak{f}(\chi),
  \]
  where $\chi$ is any algebraic Hecke character of conductor dividing $\mathfrak{f}\pi^\infty$ with infinite type $(-k, 0)$ for $k>0$, and $G(\chi)$ is the Gauss sum.
  \item[$(2)$.] For the measure $\overline{\mu_\mathfrak{f}}$ on $\Gal(K(\pi^\infty)/K)$ induced by $\mu_{\mathfrak{f}}$ via the natural projection $\mathcal{G}_\mathfrak{f} \twoheadrightarrow \Gal(K(\pi^\infty)/K)$, we have
  \[
  \overline{\mu_\mathfrak{f}} = \prod_{\mathfrak{l} \mid \mathfrak{f}} (1 - \mathrm{Frob}_{\mathfrak{l}}^{-1}) \cdot \mu_1,
  \]
  where $\mu_1 \in \mathrm{Frac}(\Lambda_{\Gal(K(\pi^\infty)/K)} \otimes_{\mathbb{Z}_p} W)$ does not depend on $\mathfrak{f}$. Moreover, for any $\sigma \in \Gal(K(\pi^\infty)/K(\pi))$, we have $(1 - \sigma) \mu_1 \in \Lambda_{\Gal(K(\pi^\infty)/K)} \otimes_{\mathbb{Z}_p} W$.
\end{enumerate}
\end{theorem}
\begin{remark}
\begin{enumerate}
\item $\mu_\mathfrak{f}$ does not depend on the choice of $\omega$. Indeed, if we replace $\omega$ by $u\omega$, where $u \in \O_v^\times$, $(\Omega, \Omega_p)$ changes to $(u\Omega, u\Omega_p)$.
\item Let $\mathfrak{g}$ be the integral ideal such that $\mathfrak{g}|\mathfrak{f}$. Then, we have 
$$L_\mathfrak{f}(\chi)=\prod_{v|\mathfrak{f}, v\not|\mathfrak{g}}(1-\chi_{\mathrm{gal}}(\mathrm{Frob}_v))L_\mathfrak{g}(\chi).$$
From this relation, we obtain the second half of the claim.
\end{enumerate}
\end{remark}

Now, let $\pi: \mathrm{Frac}(W\hoge[\Gal(K(\pi^\infty)/K)]) \to \mathrm{Frac}(W\hoge[\Gamma])$ be the field homomorphism induced by the natural projection $\Gal(K(\pi^\infty)/K) \to \Gamma$. Then define
$
\mu := \pi(\mu_1) \in \mathrm{Frac}(W\hoge[\Gamma]).
$

\begin{lemma}[\cite{ds}, Proposition 1.12(ii)]
$\mu \not\in \Lambda_{\Gamma} \otimes_{\Z_p} W$.
\end{lemma}
Let $F$ be a finite abelian extension of $K$ such that there is an inclusion $\Delta:=\Gal(F/K)\hookrightarrow (\Z/p)^\times \times (\Z/p)^\times$. Let $K_\infty$ be the maximal $\Z_p$-extension unramified outside $\pi$, and denote the field composition $F_\infty:=FK_\infty$.

\begin{defin}\label{def:p-L}
For the above $F$, let $\mathfrak{f}$ be the prime-to-$\pi$ part of the modulus of $F$. Then, we have the restriction map $\pi_F:W\hoge[\Gal(K(\mathfrak{f}\pi^\infty)/K)]\twoheadrightarrow W\hoge[F_\infty/K]$. 
For $\mathfrak{f}\neq 1$ case, we denote $\mu_F:=\pi_F(\mu_{\mathfrak{f}})\in W\hoge[\Gal(F_\infty/K)].$

For $\mathfrak{f}=1$ case, we always fix a topological generator of $\Gal(K_\infty/K)$ denoted $\gamma_0$ and define $\mu_F:=(\gamma_0-1)\mu\in \Lambda_\Gamma\otimes_{\Z_p}W$.

Furthermore, for each $\chi\in \widehat{\Delta}=\widehat{\Gal(F/K)}$, let $\mu_F(\chi):=e_\chi \mu_F\in \Lambda_\Gamma\otimes_{\Z_p}W$, where $e_\chi:=1/|\Delta|\sum_{g\in \Delta}\chi(g)g^{-1}$ is the idempotent attached to $\chi$.
\end{defin}
For the two-variable case, we define as follows:
\begin{defin}
Let \(F/K\) be as above, and let \(K_\infty/K\) be the maximal
\(\mathbb Z_p^2\)-extension of \(K\). Put \(F_\infty:=FK_\infty\).
Let \(\mathfrak f\) be the prime-to-\(p\) part of the modulus of \(F/K\).
We denote by
\[
\pi_F:W\hoge[\operatorname{Gal}(K(\mathfrak f p^\infty)/K)]
\longrightarrow
W\hoge[\operatorname{Gal}(F_\infty/K)]
\]
the natural homomorphism. Then we define the two-variable \(p\)-adic
\(L\)-function associated with \(\pi\) by
\(
\mu_{F,\pi}:=\pi_F(\mu_{\mathfrak f\bar\pi^\infty}).
\)

Furthermore, for each $\chi\in \widehat{\Delta}$, we define $\mu_{F,\pi}(\chi):=e_\chi\mu_{F,\pi}$, where $e_\chi=1/|\Delta|\sum_{g\in \Delta}\chi(g)g^{-1}$.
\end{defin}

For the $\Z_p$ or $\Z_p^2$-extension $K_\infty/K$, we denote $F_\infty=FK_\infty$, and define $\mathfrak{X}_{F,\pi} := \Gal(M_\pi(F)/F_\infty)$ and $\mathfrak{X}_{K,\pi} := \Gal(M_\pi(K)/K_\infty)$, where $M_\pi(F)$ and $M_\pi(K)$ denote the maximal abelian pro-$p$ extensions of $F_\infty$ and $K_\infty$ unramified outside $\pi$, respectively. 
We introduce the necessary notation.

\begin{defin}
For $F_\infty/K$ as above, we denote $\mathcal{G}:=\Gal(F_\infty/K)$. Let $X$ be any $\Lambda_{\mathcal{G}}$-module, and let $\chi \in \widehat{\Delta}$ be any character. Define the ring homomorphism naturally induced by $\chi$ as $\chi: \Lambda_{\mathcal{G}} \to \Lambda_{\Gamma} \otimes_{\mathbb{Z}_p} W$. The $\Lambda_{\Gamma} \otimes_{\mathbb{Z}_p} W$-module $M^{\chi}$ is defined as
\[
M^{\chi} := \left\{ m \in M \otimes_{\mathbb{Z}_p} W\mid g \cdot m = \chi(g) m \text{ for all } g \in \Delta \right\}.\]
\end{defin}
The following theorem is Rubin's main conjecture for imaginary quadratic fields, proved by Rubin \cite{Rubin91}, \cite{Rubin1994MoreMainConjectures}.

\begin{theorem}(Rubin \cite{Rubin91}, \cite{Rubin1994MoreMainConjectures})\label{Conj}
Let $F/K$, and $F_\infty/K$ as above. Then, the following holds for all $\chi \in \widehat{\Delta}=\widehat{\Gal(F/K)}$:
\[
\operatorname{char}_{\Lambda_{\Gamma}\widehat{\otimes}_{\mathbb Z_p} W}
\bigl(\mathfrak X_{F,\pi}^{\chi}\bigr)
=
\begin{cases}
\bigl(\mu_F(\chi)\bigr),
& \text{if } K_\infty/K \text{ is the } \mathbb Z_p\text{-extension}, \\[2mm]
\bigl(\mu_{F,\pi}(\chi)\bigr),
& \text{if } K_\infty/K \text{ is the } \mathbb Z_p^2\text{-extension}.
\end{cases}
\]
\end{theorem}
\begin{remark}
Rubin proved that the characteristic element of the relevant Iwasawa module coincides with that of the quotient of the product of local units by elliptic units. By a slight modification of de Shalit’s construction \cite{ds}, this characteristic element can be identified with the $p$-adic $L$-function associated with \(\pi\). For a concrete reference, see Kataoka’s work \cite{Kataoka2022TwoVariableFitting}.

\end{remark}
By definition, the following surjection exists:
\[
\mathfrak{X}_{F,\pi} = \Gal(M_\pi(F)/F_\infty) \twoheadrightarrow  \mathfrak{X}_{K,\pi} = \Gal(M_\pi(K)/K_\infty).
\]

\begin{center}
\[
\xymatrix{
  M_\pi(K) \ar@{-}[r] & M_\pi(F) \\
  K_\infty \ar@{-}[r] \ar@{-}[u] & F_\infty  \ar@{-}[u]
}
\]
\end{center}

Since $p\nmid |\Delta|$, put
\[
e_\Delta:=\frac{1}{|\Delta|}\sum_{\delta\in\Delta}\delta.
\]
Then
\[
\mathfrak{X}_{F,\pi}^\Delta=e_\Delta\mathfrak{X}_{F,\pi}
\simeq
\mathfrak{X}_{F,\pi}/(1-e_\Delta)\mathfrak{X}_{F,\pi}.
\]
Let $M^\Delta$ be the fixed field of $(1-e_\Delta)\mathfrak{X}_{F,\pi}$ in
$M_\pi(F)$. Then
\[
\operatorname{Gal}(M^\Delta/F_\infty)\simeq \mathfrak{X}_{F,\pi}^\Delta.
\]
Moreover, we have the following:

\begin{lemma}\label{Iwasawa-mod-Imaginary}
The natural map $\mathfrak{X}_{F,\pi}^{\Delta} \to \mathfrak{X}_{K,\pi}$ is an isomorphism. In particular, we have $(\mu_F(\chi_{\mathrm{triv}}))=(\mu_K)$.
\end{lemma}

\begin{proof}
By construction, the action of $\Delta$ on
$\operatorname{Gal}(M^\Delta/F_\infty)\simeq e_\Delta\mathfrak{X}_{F,\pi}$ is trivial.
Since $F_\infty/K_\infty$ is abelian with Galois group $\Delta$, it follows that
$\operatorname{Gal}(M^\Delta/K_\infty)$ is abelian. Moreover,
$M^\Delta/K_\infty$ is unramified outside the primes above $\pi$.
The exact sequence $$1\rightarrow P^{'}:=\Gal(M^{\Delta}/F_\infty)\rightarrow \Gal(M^{\Delta}/K_\infty)\twoheadrightarrow \Gal(F_\infty/K_\infty)\simeq \Delta\rightarrow 1$$
splits because $|\Delta|$ is coprime to $p$. So, we have $\Gal(M^\Delta/K_\infty)\simeq P^{'}\times \Delta$. Let $L_\Delta$ be the field corresponding to the subgroup $\Delta\subset \Gal(M^{\Delta}/K_\infty)$. Then, $L_\Delta/K_\infty$ is pro-$p$ extension unramified outside $\pi$. Hence, $L_\Delta\subset M_\pi(K)$ and therefore, 
\[
M^\Delta\subset M_\pi(K)F_\infty.
\]
Conversely, $M_\pi(K)F_\infty/F_\infty$ is obtained by base change from
$K_\infty$, and its Galois group over $F_\infty$ is fixed by $\Delta$. Hence
\[
M_\pi(K)F_\infty\subset M^\Delta.
\]

Thus $M^\Delta=M_\pi(K)F_\infty$, and 
\[
\operatorname{Gal}(M^\Delta/F_\infty)\simeq\mathfrak{X}_{K,\pi}.
\]
This proves the lemma.
\end{proof}

From the Iwasawa Main Conjecture for imaginary quadratic fields, we obtain the following corollary:

\begin{corollary}\label{Iwasawa-module}
Let $K$ be an imaginary quadratic field such that $\mathrm{Cl}_K\otimes \Z_p$ is cyclic, and let $K_\infty/K$ be the $\Z_p$-extension unramified outside $\pi$. Then, we have 
\[\mathfrak{X}_{K,\pi} \otimes_{\Z_p} W\simeq \Lambda_{\Gamma} \otimes_{\Z_p} W / (\mu_K).
\]
\end{corollary}

\begin{proof}
By Proposition~\ref{Iwasawa-module-top}, the \(\Lambda_{\mathbb Z}\)-module
\(\mathfrak{X}_{K,\pi}\) is cyclic.  Hence
\(
\mathfrak{X}_{K,\pi}\simeq \Lambda_{\mathbb Z}/I
\)
for some ideal \(I\subset\Lambda_{\mathbb Z}\).

By Greenberg's theorem recalled in Proposition~\ref{General-Greenberg}, this Iwasawa module
contains no non-trivial pseudo-null submodule.  On the other hand, Theorem \ref{Conj} gives
\[
\operatorname{char}_{\Lambda_W}
\bigl(\mathfrak{X}_{K,\pi}\widehat{\otimes}_{\mathbb Z_p}W\bigr)
=
(\mu_K).
\]
Let \(f_K\in\Lambda_{\mathbb Z}\) be a generator of the divisorial closure
of \(I\).  Then
\(
(f_K)\Lambda_W=(\mu_K).
\)
We claim that \(I=(f_K)\).  Indeed, for every height-one prime
\(\mathfrak p\subset\Lambda_{\mathbb Z}\), the equality of characteristic
ideals implies
\(
I_{\mathfrak p}=(f_K)_{\mathfrak p}.
\)
Hence \((f_K)/I\) is supported in codimension at least two, and therefore is
pseudo-null.  But \((f_K)/I\) is a submodule of
\(\Lambda_{\mathbb Z}/I\simeq\mathfrak{X}_{K,\pi}\), which has no non-trivial
pseudo-null submodule.  Therefore \((f_K)/I=0\), so \(I=(f_K)\).

After tensoring with \(W\), we obtain
\[
\mathfrak{X}_{K,\pi}\widehat{\otimes}_{\mathbb Z_p}W
\simeq
\Lambda_W/(f_K)
\simeq
\Lambda_W/(\mu_K).
\]
\end{proof}
Finally, we will study the $\Z_p^2$-extension $K_\infty/K$, and the Iwasawa module of $F_\infty=K_\infty F$.

Consider the canonical exact sequence of
\(\mathbb Z_p\hoge[\Delta\times\Gamma]\)-modules
\[
0
\longrightarrow
Y_{\bar{\pi}}
\longrightarrow
\mathfrak{X}_{F,p}
\longrightarrow
\mathfrak{X}_{F,\pi}
\longrightarrow
0.
\]
In the same manner, we define  $Y_{\pi}$.
Now, we introduce a Fitting ideal. 
\begin{defin}
Let $R$ be a commutative ring and let $M$ be a finitely presented $R$-module. Choose a finite presentation
\[
R^m\stackrel{\Psi_M}{\longrightarrow}R^n\twoheadrightarrow M.
\]
Define $\mathrm{Fitt}_{i}(M)\subset R$ to be the ideal of $R$ generated by the determinants of all $(n-i)\times (n-i)$ submatrices of $\Psi_M$. This does not depend on the choice of the presentation of $M$.
\end{defin}

\begin{lemma}[Northcott~\cite{Northcott76}, or Traldi~\cite{Traldi82} Lemma 2.1(d) ]\label{lemma:Northcot}
Let $R$ be a ring, and 
 \(0\to A\to B\to C\to0\) be an exact sequence of 
finitely presented \(R\)-modules. Then,  the following inclusion holds for $k\ge 0$;
\[
\sum_{i+j=k}\operatorname{Fitt}_i(A)\operatorname{Fitt}_j(C)
\subset
\operatorname{Fitt}_k(B).
\]

\end{lemma}
From this lemma, we obtain the following corollary.
\begin{corollary}\label{cor:p-iwasawa}
For the above $\mathfrak{X}_{F,p}\widehat{\otimes}_{\Z_p}W$ and $\chi\in \widehat{\Delta}$, we have 
$$\mathrm{Fitt}_{1}(\mathfrak{X}_{F,p}^{\chi}\widehat{\otimes}_{\Z_p}W)\supset \mu_{F, \pi}(\chi)\mathrm{Fitt}_{1}(Y_{\bar\pi}^{\chi}\widehat{\otimes}_{\Z_p}W)+\mu_{F,\bar\pi}(\chi)\mathrm{Fitt}_{1}(Y_{\pi}^{\chi}\widehat{\otimes}_{\Z_p}W).$$
\end{corollary}
The cyclicity of $Y_{\bar\pi}^{\chi}$ as a $\Lambda_{\Gamma}$-module is well studied in \cite{BCGKPPST} in the case when $\zeta_p\in F$. We can apply their claim even if $\zeta_p\not\in F$ with a slight modification.

\begin{prop}\label{prop:p-iwasawa-2}

Let
$
        \Delta=\operatorname{Gal}(F/K)
$
and let \(\chi\) be a character of \(\Delta\).  Fix a prime
\(\mathfrak P\) of \(F\) above \(\bar{\pi}\), and denote by
\(
        \Delta_{\mathfrak P}
        \subset \Delta
\)
its decomposition subgroup.  Let
\(
        Y_{\bar\pi}^{\chi}
\)
be the \(\chi\)-component of $Y_{\bar\pi}$, viewed as a module over
\(
        \Lambda=\Z_p\hoge[\Gamma].
\)

Assume either
\(
        F_{\mathfrak P}\not\supset K_{\bar{\pi}}(\zeta_p),
\)
or
\[
        F_{\mathfrak P}\supset K_{\bar{\pi}}(\zeta_p)
        \quad\text{and}\quad
        \chi|_{\Delta_{\mathfrak P}}\neq \omega_{\mathfrak P},
\]
where \(\omega_{\mathfrak P}\) denotes the local mod-\(p\) cyclotomic
character on \(\Delta_{\mathfrak P}\).  Then
\(
        Y_{\bar \pi}^{\chi}
        \simeq
        \Lambda
\)
as \(\Lambda\)-modules.  
\end{prop}
\begin{proof}
Put
\[
        \widetilde F=F(\zeta_p),
        \qquad
        \widetilde F_\infty=\widetilde F K_\infty,
        \qquad
        \widetilde\Delta=\operatorname{Gal}(\widetilde F/K),
\]
and let
\(
        H=\operatorname{Gal}(\widetilde F/F).
\)
Since \([\widetilde F:F]\) is prime to \(p\), the extension
\(\widetilde F_\infty/F_\infty\) is finite of degree prime to \(p\), and
\(
        \operatorname{Gal}(\widetilde F_\infty/\widetilde F)
        \simeq
        \Gamma.
\)
Let \(\widetilde\chi:\widetilde{\Delta}\rightarrow \Delta \stackrel{\chi}{\rightarrow }\C^\times\) be the inflation of \(\chi\) to
\(\widetilde\Delta\).

Choose a prime \(\widetilde{\mathfrak P}\) of \(\widetilde F\) above
\(\mathfrak P\), and let
\(
        \widetilde\Delta_{\widetilde{\mathfrak P}}
        \subset \widetilde\Delta
\)
be the decomposition subgroup at \(\widetilde{\mathfrak P}\).  We have an
exact sequence
\[
        1
        \longrightarrow
        H_{\widetilde{\mathfrak P}}
        \longrightarrow
        \widetilde\Delta_{\widetilde{\mathfrak P}}
        \longrightarrow
        \Delta_{\mathfrak P}
        \longrightarrow
        1,
\]
where
\(
        H_{\widetilde{\mathfrak P}}
        \simeq
        \operatorname{Gal}
        \bigl(F_{\mathfrak P}(\zeta_p)/F_{\mathfrak P}\bigr).
\)

Step 1. we prove that $Y_{\bar\pi}\simeq \prod_{v|\bar\pi}
I_v$ holds. Here, $I_v$ is an inertia group of a maximal abelian pro-$p$ extension of $F_{\infty,v}$.

Let $X_{\mathrm{ur}}$ be the unramified Iwasawa module of $F_\infty$. 
For a finite extension \(L/F\) contained in \(F_\infty\), put
\[
\mathcal U_p(L)
:=
\prod_{v\mid p}\widehat{\mathcal O_{L_v}^{\times}},
\qquad
\mathcal U_{\pi}(L)
:=
\prod_{v\mid \pi}\widehat{\mathcal O_{L_v}^{\times}},
\qquad
\mathcal U_{\bar\pi}(L)
:=
\prod_{v\mid \bar\pi}\widehat{\mathcal O_{L_v}^{\times}}.
\]
We then define
\[
\mathcal U_p
:=
\varprojlim_{L\subset F_\infty}\mathcal U_p(L),
\qquad
\mathcal U_{\pi}
:=
\varprojlim_{L\subset F_\infty}\mathcal U_{\pi}(L),
\qquad
\mathcal U_{\bar\pi}
:=
\varprojlim_{L\subset F_\infty}\mathcal U_{\bar\pi}(L),
\]
where the transition maps are induced by the local norm maps.
Equivalently,
\(
\mathcal U_p
=
\mathcal U_{\pi}\times\mathcal U_{\bar\pi}.
\)
Similarly, define $\mathcal{E}_p$, $\mathcal{E}_\pi$ be the inverse limit of closure of the image $\O_L^\times \rightarrow \mathcal{U}_p(L)$, $\O_L^\times \rightarrow \mathcal{U}_\pi(L)$ for finite subfields $L\subset F_\infty $ respectively. 
Class field theory gives a commutative diagram with exact rows
\[
\begin{array}{ccccccccc}
0 &\longrightarrow&
\mathcal U_p/\overline{\mathcal E}_p
&\longrightarrow&
\mathfrak X_p
&\longrightarrow&
X_{\mathrm{ur}}
&\longrightarrow& 0
\\
&& \downarrow && \downarrow && \Vert &&
\\
0 &\longrightarrow&
\mathcal U_{\pi}/\overline{\mathcal E}_{\pi}
&\longrightarrow&
\mathfrak X_{\pi}
&\longrightarrow&
X_{\mathrm{ur}}
&\longrightarrow& 0 .
\end{array}
\]
Put \(
Y_{\bar\pi}
\simeq
\ker\left(
\mathcal U_p/\overline{\mathcal E}_p
\longrightarrow
\mathcal U_{\pi}/\overline{\mathcal E}_{\pi}
\right)
\). By the snake lemma, we obtain a $\Lambda$-module exact sequence
$$0\rightarrow Y_{\bar\pi}\rightarrow \mathfrak{X}_p\rightarrow \mathfrak{X}_\pi\rightarrow 0$$

The map
\[
\prod_{v\mid\bar\pi} I_v
\simeq
\mathcal U_{\bar\pi}
\longrightarrow
Y_{\bar\pi},
\qquad
x\longmapsto [(1,x)]
\]
is surjective.  Its kernel is
\(
\overline{\mathcal E}_p
\cap
\bigl(\{1\}\times \mathcal U_{\bar\pi}\bigr).
\)
Indeed, \((1,x)\) maps to zero in \(Y_{\bar\pi}\) if and only if
\((1,x)\in\overline{\mathcal E}_p\).

By the \(\Sigma\)-Leopoldt theorem for abelian extensions of imaginary
quadratic fields, applied to
\(
\Sigma=\{v\mid\pi\},
\)
the one-sided localization map
\(
\overline{\mathcal E}_p
\longrightarrow
\mathcal U_{\pi}
\)
is injective.  Therefore
\(
\overline{\mathcal E}_p
\cap
\bigl(\{1\}\times \mathcal U_{\bar\pi}\bigr)
=0.
\)
Consequently, the above surjection is an isomorphism, and hence
\(
Y_{\bar\pi}
\simeq
\mathcal U_{\bar\pi}
\simeq
\prod_{v\mid\bar\pi} I_v .
\)
We define \(\widetilde{Y}_{\bar\pi}\) in the same way as \(Y_{\bar\pi}\),
replacing \(F\) by \(\widetilde F\) and \(F_\infty\) by \(\widetilde F_\infty\).

Step 2. Application of \cite{BCGKPPST}.

 We define two modules $D_{\bar\pi}$, and $\widetilde{D_{\bar{\pi}}}$ as follows;
$$D_{\bar{\pi}}:=\prod_{v|\bar\pi}\plim[F_\infty\supset L\supset F]D_{v}(L), \ \widetilde{D_{\bar{\pi}}}:=\prod_{v|\bar\pi}\plim[\widetilde{F_\infty}\supset L\supset \widetilde{F}]D_{v}(L).$$
Here, $v\mid \bar\pi$ is the finite place of $F_\infty$ or $\widetilde{F_\infty}$ and $D_{v}(L)$ is the decomposition group of $L/K$ at $v\mid \bar\pi$. Note that the number of places of $F_\infty$ and $\widetilde{F_\infty}$ over $\bar\pi$ is finite.
Bleher--Chinburg--Greenberg--Kakde--
Pappas--Sharifi--Taylor studied the $\Lambda_{\Gal(\widetilde{F_\infty}/K)}$-module structure of $\widetilde{D}_{\bar\pi}$ which is canonically isomorphic to $\widetilde{Y_{\bar\pi}}$. See \cite[Lemma 5.2.2]{BCGKPPST} and Step 1 above for this canonical isomorphism. In the following, for $\Lambda$-modules $M$, we use $M^{*}:=\mathrm{Hom}_\Lambda(M, \Lambda)$. Remark that we have a canonical $\Lambda$-module morphism $M\rightarrow M^{**}:m\mapsto \Tilde{m}$ where $\Tilde{m}(f):=f(m)$ for all $f\in M^{*}$.

In \cite[Lemma 4.3.3]{BCGKPPST}, it is proved that if
\(        \omega\widetilde\chi^{-1}
        |_{\widetilde\Delta_{\widetilde{\mathfrak P}}}
        \neq 1,
\)
then the canonical morphism
\[
        \widetilde{D}_{\bar{\pi}}^{\widetilde\chi}\simeq \widetilde{Y}_{\bar{\pi}}^{\widetilde\chi}
        \longrightarrow
        \left(
        (\widetilde{Y}_{\bar\pi})^{\widetilde\chi}
        \right)^{**}\simeq \left(\widetilde{D}_{\bar{\pi}}^{\widetilde\chi}\right)^{**}
\]
is an isomorphism.  This module has \(\Lambda\)-rank one and
\(\Lambda\) is a regular local UFD; its double dual is free of rank one.  we prove this in a much more general setting.

Step 3. Application of ring theory.
Let $R$ be a regular local UFD ring, and let $M$ be an $R$-module such that $M\otimes_R \mathrm{Frac}(R)$ is one-dimensional over $\mathrm{Frac}(R)$. Then we prove that $M^{**}$ is a free $R$-module of rank one. First, $N:=M^{**}$ is a reflexive $R$-module, i.e. $N^{**}\simeq N$. See, \cite[More on Algebra, Definition~15.24.9 and Lemma~15.24.8,
Tags~0AV4 and~0AV3]{stacks-project}.
By Serre's criterion for reflexive modules, we have $N=\bigcap_{\mathfrak{p}:\text{ht $1$ prime}}N_{(\mathfrak{p})}$, and $N$ is a torsion-free $R$-module. $N_{(\mathfrak{p})}$ is a module over the discrete valuation ring $R_{(\mathfrak{p})}$, which is torsion-free over $R_{(\mathfrak{p})}$, $N_{(\mathfrak{p})}\otimes_{R_{(\mathfrak{p})}}\mathrm{Frac}(R)=\mathrm{Frac}(R)$, so we obtain $N_{(\mathfrak{p})}\simeq R_{(\mathfrak{p})}$.

After fixing a non-zero element of $N$, we obtain $N=\bigcap_{\mathfrak{p}:\text{ht $1$ prime}}N_{\mathfrak{p}}\hookrightarrow \mathrm{Frac}(R)$ as an $R$-module. Hence $N$ is a fractional ideal of $R$.  Furthermore, $R$ is a UFD, so we obtain that $N$ is generated by one element. Hence, we obtain $N\simeq R$.
In particular, for $R=\Lambda$, and $M=\widetilde{Y}_{\bar\pi}^{\widetilde{\chi}}$, we obtain
\(
\widetilde{Y}_{\bar\pi}^{\widetilde{\chi}}        \simeq
        \Lambda.
\)

Here, we will check $\widetilde{Y}_{\bar\pi}^{\widetilde{\chi}}\otimes_{R}\mathrm{Frac}(R)\simeq \mathrm{Frac}{R}$. 
Put
\(Q=\operatorname{Frac}(\Lambda)\). Fix a prime \(\mathfrak P\) of \(F\)
above \(\bar\pi\), and let \(\Delta_{\mathfrak P}\) be its decomposition subgroup.
Since the primes of \(F\) above \(\bar\pi\) form a transitive \(\Delta\)-set, the
local description in Step 1 gives
\[
Y_{\bar\pi}\simeq
\operatorname{Ind}_{\Delta_{\mathfrak P}}^\Delta \widehat{\O_{F,\mathfrak{P}}^\times}.
\]
After tensoring with \(Q\), local class field theory and the normal basis
theorem give
\[
Y_{\bar\pi}\otimes_{\Lambda}Q
\simeq Q[\Delta_{\mathfrak P}].
\]
Hence
\[
Y_{\bar\pi}\otimes_{\Lambda}Q
\simeq
\operatorname{Ind}_{\Delta_{\mathfrak P}}^\Delta Q[\Delta_{\mathfrak P}]
\simeq Q[\Delta].
\]
Therefore
\[
Y_{\bar\pi}^{\chi}\otimes_{\Lambda}Q
\simeq e_\chi Q[\Delta]\simeq Q.
\]
In particular, \(Y_{\bar\pi}^{\chi}\) has \(\Lambda\)-rank one.

Step 4. Descent to $Y_{\bar\pi}$. 
 
Let
\(
q:\widetilde{\Delta}\twoheadrightarrow \Delta
\)
be the natural quotient map, and put
\(
\widetilde{\chi}:=\chi\circ q.
\)
Then, define
\[
e_{\widetilde{\chi}}
=
\frac{1}{|\widetilde{\Delta}|}
\sum_{\tau\in\widetilde{\Delta}}
\widetilde{\chi}(\tau)\tau^{-1},\ e_H=\frac{1}{|H|}\sum_{h\in H} h .
\]
Since \(\widetilde{\chi}\) is trivial on \(H\), this idempotent $e_{\widetilde{\chi}}$ lies in the
\(e_H\)-part.
Let \(L/F\) be a finite extension contained in \(F_\infty/F\), and put
\(
\widetilde L:=L\widetilde F.
\)
For the primes above \(\bar\pi\), define
\[
U_{\bar\pi}(L):=
\prod_{w\mid\bar\pi}
\widehat{\mathcal O_{L_w}^{\times}},
\qquad
U_{\bar\pi}(\widetilde L):=
\prod_{\widetilde w\mid\bar\pi}
\widehat{\mathcal O_{\widetilde L_{\widetilde w}}^{\times}},
\]
where the hat denotes \(p\)-completion.

There are canonical maps
\(
\iota_{\widetilde L/L}:
U_{\bar\pi}(L)\longrightarrow U_{\bar\pi}(\widetilde L)
\) 
and
\(
N_{\widetilde L/L}:
U_{\bar\pi}(\widetilde L)\longrightarrow U_{\bar\pi}(L)
\)
defined as follows.  The map \(\iota_{\widetilde L/L}\) is induced by the
local inclusions
\[
L_w\hookrightarrow \widetilde L_{\widetilde w}
\qquad
(\widetilde w\mid w).
\]
The norm map is given by
\[
N_{\widetilde L/L}\bigl((x_{\widetilde w})_{\widetilde w}\bigr)
=
\left(
\prod_{\widetilde w\mid w}
N_{\widetilde L_{\widetilde w}/L_w}(x_{\widetilde w})
\right)_{w}.
\]
Then
\(
N_{\widetilde L/L}\circ \iota_{\widetilde L/L}
=
[\widetilde L:L]
\)
on \(U_{\bar\pi}(L)\), and
\(
\iota_{\widetilde L/L}\circ N_{\widetilde L/L}
=
\sum_{h\in H}h
\)
on \(U_{\bar\pi}(\widetilde L)\).  Since \(|H|=[\widetilde L:L]\) is prime
to \(p\), it is a unit in \(\mathbb Z_p\).  Hence
\(
\iota_{\widetilde L/L}:
U_{\bar\pi}(L)
\xrightarrow{\ \sim\ }
e_H U_{\bar\pi}(\widetilde L).
\)

Moreover, the map \(\iota_{\widetilde L/L}\) is compatible with the quotient
maps
\(
\operatorname{Gal}(\widetilde L/K)\twoheadrightarrow
\operatorname{Gal}(L/K)
\)
and
\(
\widetilde\Delta\twoheadrightarrow\Delta.
\)
Therefore, it induces an isomorphism 
\(
e_{\chi}\bigl(U_{\bar\pi}(L)\bigr)
\xrightarrow{\ \sim\ }
e_{\widetilde\chi}
\bigl(U_{\bar\pi}(\widetilde L)\bigr).
\)
Passing to the inverse limit over all finite extensions \(L/F\) contained in
\(F_\infty/F\), we obtain
\(
U_{\bar\pi}(F_\infty)^\chi
\simeq
U_{\bar\pi}(\widetilde F_\infty)^{\widetilde\chi}.
\)

By Step 1, we have canonical class-field-theoretic identifications
\[
Y_{\bar\pi}\simeq U_{\bar\pi}(F_\infty),
\qquad
\widetilde Y_{\bar\pi}\simeq U_{\bar\pi}(\widetilde F_\infty).
\]
Therefore
\(
Y_{\bar\pi}^{\chi}
\simeq
\widetilde Y_{\bar\pi}^{\widetilde\chi}.
\)
By Step 3, the right-hand side is isomorphic to \(\Lambda\). Hence
\(
Y_{\bar\pi}^{\chi}\simeq\Lambda.
\)
\end{proof}

\section{Relations among Topological Generators}
In this section, we study the structure of the Galois group $G^\pi$ over $K$ for the maximal pro-$p$ extension of $K$ unramified outside $\pi$ and prove an imaginary quadratic analogue of Komatsu's theorem in \cite{Kom} for real quadratic fields. Here again, we use $W=W(\bar{\F}_p)$.
\begin{theorem}\label{Main 1}
Let $K$ be an imaginary quadratic field. Let $\pi$ be a prime ideal above an odd prime $p$ which splits in $K/\Q$. Let $K^\pi$ be the maximal pro-$p$ extension of $K$ unramified outside $\pi$, and let $K_\infty$ be the maximal $\Z_p$-subextension of $K^\pi/K$. Define $\Gamma := \Gal(K_\infty/K)$. Let $\mu_K \in W\hoge[\Gamma]$ be the $p$-adic $L$-function of $K$.

Let $K_p^{(1)}$ denote the maximal $p$-subextension of the Hilbert class field of $K$. Assume that the following conditions hold for $(K, p)$:
\begin{enumerate}
    \item[$1$.] The $p$-part of the ideal class group $\mathrm{Cl}_K \otimes \Z_p$ is cyclic.
    \item[$2$.] $K_p^{(1)}\not\subset K_\infty$.
\end{enumerate}
Then there exists a pro-$p$ free group $F = \langle \gamma, g \rangle$ of rank 2 and a non-trivial element $r \in F$ such that $G^\pi = \Gal(K^\pi/K) \simeq F / \overline{\langle r \rangle}$. Here, $\gamma$ is a topological generator of $\Gamma$, and $\overline{\langle r \rangle}$ is the minimal normal closed subgroup containing $r$. Furthermore, if $N$ is the smallest normal subgroup of $F$ containing $g$, then $N / [N, N] \simeq\Z_p\hoge[\Gamma]$. Under the isomorphism $N/[N,N]\simeq\Z_p\hoge[\Gamma]$ sending the image of $g$ to $1$, for some $u \in W\hoge[\Gamma]^\times$, we have
\[
r \equiv u \mu_K \mod [N, N] \widehat{\otimes}_{\Z_p} W.
\]
In particular, $r\mod [N,N]\in N/[N,N]\simeq \Lambda_\Gamma$ is coprime to $p$ in $\Lambda_\Gamma$.
\end{theorem}

We first introduce the following lemma on pro-$p$ groups, proved by Komatsu~\cite{Kom}.

\begin{lemma}[\cite{Kom}]\label{two-free}
Let $F = \langle x, y \rangle$ be a free pro-$p$ group of rank 2. Assume that the closed subgroup topologically generated by $y$ is isomorphic to $\Z_p$. Let $N$ be the smallest closed normal subgroup of $F$ containing $x$. Then there is an isomorphism
\[
\Phi: N / [N, N] \simeq \Z_p\hoge[\Gamma],
\]
induced by $\Phi(x [N, N]) = 1$.
\end{lemma}

\begin{proof}[Proof of Theorem~\ref{Main 1}]
The topological generating system is characterized by Proposition~\ref{Iwasawa-module-top}. 
We recall the relation-rank estimate in the mixed ramification case.
Let $S=\{\pi\}$ and
\[
\mathcal V_S(K):=
\left\{
a\in K^\times/K^{\times p}\ \middle|\ 
\begin{array}{l}
a\in K_{\mathfrak q}^{\times p}\ \text{for every } \mathfrak q\in S,\\
v_{\mathfrak q}(a)\equiv 0 \mod p\ \text{for every } \mathfrak q\notin S
\end{array}
\right\},
\]
 and let 
\(
\mathcal B_S(K):=\operatorname{Hom}_{\mathbb F_p}
(\mathcal V_S(K),\mathbb F_p).
\)
By the proof of Vogel's mixed-case theorem
\cite[Proposition~2.2 and its proof]{Vogel2007Mixed}, one has
\[
h^2(G^\pi)
\le
\sum_{\mathfrak q\in S}\delta_{\mathfrak q}
+
\dim_{\mathbb F_p}\mathcal B_S(K),
\]
where
\(
h^2(G^\pi)
:=
\dim_{\mathbb F_p}H^2(G^\pi,\mathbb F_p)
\)
and
\[
\delta_{\mathfrak q}
=
\begin{cases}
1, & \text{if } \mu_p\subset K_{\mathfrak q},\\
0, & \text{otherwise}.
\end{cases}
\]

Since the odd prime \(p\) splits in \(K\), we have 
\(
K_\pi\simeq \mathbb Q_p,
\)
 which does not contain \(\mu_p\).  Hence
\(
\delta_\pi=0.
\)
Therefore
\[
h^2(G_S(K)(p))
\le
\dim_{\mathbb F_p}\mathcal B_S(K).
\]

We now bound \(\mathcal B_S(K)\).  Since \(S=\{\pi\}\), the defining local
condition at \(\pi\) gives an inclusion
\[
\mathcal V_S(K)\subset \mathcal V_{\emptyset}(K).
\]
There is a standard exact sequence
\[
0
\longrightarrow
\mathcal O_K^\times/\mathcal O_K^{\times p}
\longrightarrow
\mathcal V_{\emptyset}(K)
\longrightarrow
\operatorname{Cl}_K[p]
\longrightarrow
0.
\]
Because \(K\) is imaginary quadratic and \(p\) is an odd split prime, we have
\(
p\nmid |\mathcal O_K^\times|.
\)
Thus
\(
\mathcal O_K^\times/\mathcal O_K^{\times p}=0,
\)
and consequently
\(
\mathcal V_{\emptyset}(K)\simeq \operatorname{Cl}_K[p].
\)
Since \(A_K\) is cyclic, \(\operatorname{Cl}_K[p]\) has
\(\mathbb F_p\)-dimension at most \(1\).  Hence
\[
\dim_{\mathbb F_p}\mathcal B_S(K)
=
\dim_{\mathbb F_p}\mathcal V_S(K)
\le 1.
\]
It follows that
\(
h^2(G^\pi)\le 1.
\)

It remains to exclude the possibility \(h^2(G^\pi)=0\), i.e. \(G^\pi\) is a free pro-\(p\) group.  Assume that $G^\pi$ is a free pro-$p$ group, then
\(
(G^\pi)^{\mathrm{ab}}
\)
would be a free \(\mathbb Z_p\)-module, and hence would have no non-trivial
\(p\)-power torsion.

On the other hand, class field theory gives
\(
(G^\pi)^{\mathrm{ab}}_{\mathrm{tors}}\simeq B_K,
\)
where \(B_K\) is the subgroup of \(A_K\) defined as the image of
\((G^\pi)^{\mathrm{ab}}_{\mathrm{tors}}\) in \(A_K\).  By Proposition~\ref{Iwasawa-module-top},
\[
|B_K|
=
\frac{|A_K|}
{[K_\infty\cap K^{(1)}_p:K]}.
\]
Since \(A_K\) is cyclic and \(K^{(1)}_p\not\subset K_\infty\), we have
\[
[K_\infty\cap K^{(1)}_p:K]<[K^{(1)}_p:K]=|A_K|.
\]
Thus
\(
B_K\ne 0.
\)
Therefore \((G^\pi)^{\mathrm{ab}}\) has non-trivial \(p\)-power torsion,
so \(G^\pi\) cannot be free.  Hence, the relation rank of \(G^\pi\) is equal to \(1\).

First, we show that $r \in N$. By definition, $F / N \simeq \Gamma \simeq G / N$. Since $G = F / \overline{\langle r \rangle}$, if $r \not\in N$, then the natural surjection $F / N \twoheadrightarrow G / N \simeq \Gamma$ would not be an isomorphism, leading to a contradiction.

Next, we show that $r \not\in [N, N]$. Suppose $r \in [N, N]$. Let $\overline{N} := N \mod \overline{\langle r \rangle} \subset G$. Then
\[
\Z_p\hoge[\Gamma] \simeq N / [N, N] \simeq \overline{N} / [\overline{N}, \overline{N}].
\]
Here, the first isomorphism follows from Lemma~\ref{two-free}. On the other hand, $\overline{N} / [\overline{N}, \overline{N}]$ is isomorphic to the Iwasawa module $\mathfrak{X}_{K,\pi}$ defined in Corollary~\ref{Iwasawa-module}. This contradicts the fact that $\mathfrak{X}_{K,\pi}$ is a torsion $\Z_p\hoge[\Gamma]$-module, whereas $\Z_p\hoge[\Gamma]$ is a free $\Z_p\hoge[\Gamma]$-module. Hence, $r \not\in [N, N]$.

Now, let $R \in N / [N, N] \simeq \Z_p\hoge[\Gamma]$ correspond to $r \mod [N, N]$. By Proposition~\ref{Iwasawa-module-top}, $\mathfrak{X}_{K,\pi}$ is a cyclic $\Z_p\hoge[\Gamma]$-module $\Z_p\hoge[\Gamma] \cdot g$, isomorphic to $\Z_p\hoge[\Gamma] / (R)$. By Corollary~\ref{Iwasawa-module}, there exists an isomorphism
\[
W\hoge[\Gamma] / (R) \stackrel{}{\simeq} W\hoge[\Gamma] / (\mu_K).
\]
This isomorphism of cyclic $W\hoge[\Gamma]$-modules implies equality of annihilator ideals,
\(
(R)=(\mu_K)
\)
in $W\hoge[\Gamma]$. Hence $R = u \mu_K$ for some $u \in W\hoge[\Gamma]^\times$. We obtain
\[
r \equiv u \mu_K \cdot g \mod [N, N] \widehat{\otimes}_{\Z_p} W.
\]
By the \(\mu=0\) theorem due to Gillard\cite{Gillard1991} and Schneps\cite{Schneps1987}, with the correction to Gillard's argument supplied by Lamplugh\cite{Lamplugh2015}, the element \(\mu_K\) is not divisible by \(p\) in \(W\hoge[\Gamma]\). Thus $r \mod [N,N]$ is also prime to $p$ inside $N/[N,N]=\Z_p\hoge[\Gamma]\cdot g\simeq \Z_p\hoge[\Gamma]$.
\end{proof}
\begin{remark}
We emphasize that the integrality statement for the unit $u$ in the theorem is special to the present situation and is not automatic in general. Indeed, there may exist an element $f\in W(\overline{\F}_p)[\Gamma]$ such that $fu \not\in \Z_p\hoge[\Gamma]$ for all $u\in W(\overline{\F}_p)\hoge[\Gamma]$. For example, if $f=p+a(\gamma-1)$ with $a\in W(\overline{\F}_p)\setminus\Z_p$, no such unit $u$ can make $uf\in \Z_p\hoge[\Gamma]$.
\end{remark}
\section{Application to the Universal Deformation Rings}
In this section, we describe the relevant quotient deformation rings for mod $p$ Galois representations over imaginary quadratic fields. Let $S$ be either $S_\pi=\{\pi\}$ or $S_p=\{\pi,\overline{\pi}\}$. Let $K^{S}$ be the maximal algebraic extension of $K$ inside $\overline{K}$ unramified outside $S$.
We fix a mod $p$ Galois representation
$$\overline{\rho}:\Gal(K^{S}/K)\rightarrow \mathrm{GL}_2(\F_p)$$
whose image is contained in the upper triangular subgroup of $\mathrm{GL}_2(\F_p)$. We shall use the following elementary observation.
\begin{lemma}
For any subgroup $H<B(\F_p)=\{\begin{pmatrix}
* & *\\
0 & *\end{pmatrix}\}\subset \mathrm{GL}_2(\F_p)$, there exists a unipotent normal subgroup $V\triangleleft H$ such that $H/V$ is diagonal. Moreover, $H\simeq V\rtimes H/V$.
\end{lemma}
\begin{proof}
Take $V=H\cap U$, where $U$ is the maximal unipotent subgroup of $B(\F_p)$. Then $V$ is normal in $H$, because $U$ is normal in $B(\F_p)$. The quotient $H/V$ has order prime to $p$. Hence, Schur--Zassenhaus gives a splitting, and the lemma follows.
\end{proof}
Applying this lemma to $\mathrm{Im}(\overline{\rho})$, we obtain
$\mathrm{Im}(\overline{\rho})\simeq V\rtimes A$, where $V$ is unipotent and $A$ is an abelian group of order prime to $p$. We use this notation throughout this section.

Let $L$ be the finite Galois extension of $K$ corresponding to $\ker(\overline{\rho})$ and let $F$ be the subfield of $L$ corresponding to $V$. Hence, $\Delta:=\Gal(F/K)\simeq A$ holds through $\overline{\rho}$. Then $\Delta$ is a subgroup of $(\Z/p\Z)^\times\times (\Z/p\Z)^\times$.
\begin{defin}\label{def: chi_0}
For the Galois representation $\overline{\rho}$ as above, we define characters $\chi_1,\chi_2:\Delta\rightarrow \F_p^\times$ by $\overline{\rho}|_{\Delta}=\begin{pmatrix}
\chi_1 & 0\\
0 & \chi_2 \end{pmatrix}
:\Delta \rightarrow \mathrm{GL}_2(\F_p)$. We put $\chi_0:=\chi_1\chi_2^{-1}:\Delta\rightarrow \F_p^\times $.
\end{defin}
We assume now that $S=\{\pi\}$ or $=\{p\}$, and $\bar\rho$ has its image non-commutative.  Let $F_\infty$ be the maximal $\Z_p$-extension of $F$ which is unramified outside $\pi$ or $p$ and abelian over $K$, and let $M_{F,\pi}$, $M_{F,p}$ be the maximal abelian pro-$p$ extension of $F_\infty$ unramified outside the primes above $\pi$, $p$ respectively. Put $\Gamma_F:=\Gal(F_\infty/F)$.

\[
\xymatrix{
   & &M_{F,\pi} \ar@{-}[d] \\
   & F_{\infty} \ar@{-}^{\mathfrak{X}_{F,\pi}}[ur]\ar@{-}_{\Gamma_F}[d] &L \\
   & F \ar@{-}^{\Delta}[d]\ar@{-}[ur] \\
   & K
}
\]

We note that $L\subset M_{F,\pi}$ or $L\subset M_{F,p}$ according to $S=\{\pi\}$ or $=\{p\}$ because $L/F$ is an abelian extension and unramified outside all the primes of $F$ above $\pi$ or $p$.
We define $G_{F,\pi}:=\Gal(M_{F,\pi}/K)$, and $P_{F,\pi}:=\Gal(M_{F,\pi}/F)$. By Schur--Zassenhaus, we have $G_{F,\pi}\simeq P_F \rtimes \Delta$. 
Put $\mathfrak{X}_{F,\pi}:=\Gal(M_{F,\pi}/F_\infty)$.

In the same manner, we define $G_{F,p}$, and $P_{F,p}$, $\mathfrak{X}_{F,p}$.

We note that  $\Gal(F_\infty/K_\infty)\simeq \Delta$ via the canonical restriction map. So the action of $\Delta$ on $P_{F,\pi}$,  $P_{F,p}$ is always understood through this lift. Moreover, if $S=\{\pi\}$, the short exact sequence $1\rightarrow \mathfrak{X}_{F,\pi} \rightarrow P_{F,\pi}\rightarrow \Gamma_F\rightarrow 1$ splits because $\Gamma_F$ is a free pro-$p$ group. Hence $P_{F,\pi} \simeq \mathfrak{X}_{F,\pi}\rtimes \Gamma_F$ if $S=\{\pi\}$.
By our definition of $F_\infty$, which is abelian over $K$, the action of $\Delta$ on $\Gamma_F$ is trivial. Hence we obtain the following.
\begin{lemma}
With the notation above, if $S=\{\pi\}$, then $G_{F,\pi}\simeq \mathfrak{X}_{F,\pi}\rtimes (\Delta\times \Gamma_F)$.
\end{lemma}If $S=\{p\}$, we need the following lemma, which is immediate from the definition.
\begin{lemma}\label{lemma:two-var pro-p}
Let \(G\) be a pro-\(p\) group, and let \(N\subset G\) be an abelian normal subgroup. Assume that $G/N\simeq \Z_p^2=:\Gamma$. Then, $N$ becomes a $\Z_p\hoge[\Gamma]$-module by the conjugacy action of the lift. Let \(R_N\) be the set of relations among the minimal generators $n_1,\ldots, n_k$ of $N$ as $\Z_p\hoge[\Gamma]$-module, and put $c=[s,t]\in N$ where \(s,t\in G\) are lifts of topological generators of \(\Gamma\) to \(G\). Then, we have 
$G=\kakko<n_1,\ldots, n_k, s,t\ |\ R_N, [s,t]=c>.$
\end{lemma}
We now construct suitable minimal generators for $G_{F,\pi}$. By the Hochschild short exact sequence, we have $$1\rightarrow H^1(\Gamma_F,\F_p)\rightarrow H^1(P_{F,\bullet},\F_p)\rightarrow H^1(\mathfrak{X}_{F,\bullet},\F_p)^{\Gamma_F}\rightarrow H^2(\Gamma_F,\F_p),$$
for $\bullet \in \{p,\pi\}$.
Hence, when $S=\{\pi\}$, we have a short exact sequence of $\F_p[\Delta]$-modules,  $1\rightarrow (\overline{\mathfrak{X}}_{F,\pi})_{\Gamma_F}\rightarrow \overline{P}_{F,\pi}
\rightarrow \overline{\Gamma}_F \rightarrow 1 $ and when $S=\{p\}$, exact sequence $1\rightarrow \F_p \rightarrow (\overline{\mathfrak{X}}_{F,p})_{\Gamma_F}\rightarrow \overline{P}_{F,p}\rightarrow \overline{\Gamma}_F \rightarrow 1$ holds. Here, for any pro-$p$ group $C$, we denote $\overline{C}=C/C^p[C, C]$ to express its Frattini quotient. Moreover, $\F_p[\Delta]$ is semisimple and $\bar{\Gamma}_F$ has a trivial $\Delta$ action. 

We define special minimal generators of $P_F$ by using the following lemma of Boston \cite{BostonThesis}.
\begin{lemma}\label{lem:Boston}
Let $C$ be a profinite group and let $D\subset C$ be its maximal pro-$p$ normal subgroup. Put $A:=C/D$. The conjugation action of $C$ induces a representation $A\to \operatorname{Aut}(\overline D)$. For every subrepresentation $\Phi$ of $\overline D$, one may choose generators in $D$ whose images span $\Phi$ and on which $A$ acts through $\Phi$.
\end{lemma}
Then, by Lemma \ref{lem:Boston}, we obtain the following generating system of $P_F$. 
\begin{defin}\label{def:Special-gen}
We define special generators of $P_{F,\pi}$ and $P_{F,p}$, in each case $S=\{\pi\}$ and $S=\{p\}$.
\begin{enumerate}
\item In case when $S=\{\pi\}$, we denote the special generator as follows.
$$P_{F,\pi}=\langle \gamma, s_{\chi,1},\ldots,s_{\chi,s_\chi}\ |\ {\chi \in \widehat{\Delta}}\rangle$$
\begin{enumerate}
\item[$\bullet$] Let $\gamma \in P_{F,\pi}$ be the lift of the generator of $\overline{\Gamma}_F$ such that the conjugacy action of $\Delta$ is trivial. 
\item[$\bullet$] For each $\chi \in \widehat{\Delta}$, we choose elements $s_{\chi, i}\in \mathfrak{X}_{F,\pi}\subset P_{F,\pi}$ whose images form a basis of the $\chi$-isotypic part of the Frattini quotient, so that the conjugation action of $\Delta$ on $s_{\chi, i}$ is given by $\chi$.
\end{enumerate}
\item In case when $S=\{p\}$, we denote the special generator
$$P_{F,p}=\langle \gamma_1, \gamma_2, s_{\chi,1},\ldots,s_{\chi,s_\chi}\ |\ {\chi \in \widehat{\Delta}}\rangle$$
in the same manner as in $S=\{\pi\}$ case.
\end{enumerate}
\end{defin}

Now, we prove the main result for universal deformation.
\begin{proof}[Proof of Theorem \ref{Main 2}]
In the proof, put \(G:=G_{F,\bullet}\), according to whether \(\bullet=\pi\) or \(p\).
Here \(G_{F,\pi}\) and \(G_{F,p}\) are the Iwasawa-theoretic quotient groups
defined above.
First, we note that $\rho^{\mathrm{univ}}|_{\Delta}:\Delta\rightarrow \mathrm{GL}_2(R_{\bar\rho})$ is given by the Teichm{\"u}ller lift of $\bar\rho$ because $|\Delta|$ is prime to $p$.

Step 1. The description of $\rho_G(s_{\chi,i})$.

We note that for each $s_{\chi,i}$, and $g \in \Delta$, 
$gs_{\chi,i}g^{-1}=s_{\chi,i}^{\chi(g)}$ holds by definition. First, we consider the case when $\bar\rho(s_{\chi, i})=1$.
Let $M:=\rho_G(s_{\chi,i})$. First, we consider the case when $\chi$ is quadratic. If $\chi\neq 1$, choose $g\in \Delta$ such that $\chi_0(g)\neq \pm 1$, and $\chi(g)=-1$. We can always find this $g$ because $\chi_0$ is not quadratic. Then, from the equation, $\rho_G(g)M\rho_G(g)^{-1}=M^{-1}$, if we put $M=\begin{pmatrix}a&b\\
c& d\end{pmatrix}$, we have   
$$\begin{pmatrix}d& -b\\
-c& a\end{pmatrix}=\begin{pmatrix}a& \chi_0(g)b\\
\chi_0^{-1}(g)c& d\end{pmatrix}.$$
Therefore, $a=d$ and $b=c=0$, $ad=1$. So $a=d=1$ and $b=c=0$. When $\chi=1$, the same argument proves that $M$ is diagonal.

Next, we consider the case when $\chi$ is not quadratic. Fix the matrix $X:=\rho_G(g)$ where  $g\in \Delta$ satisfying $d:=\chi(g)\neq \pm 1$. Then $M\equiv 1 \mod m_R$. It is easy to check that $\det(M)=1$, $\mathrm{tr}(M)=2$ ; for completeness, we give the argument.
Indeed, from
\(
XMX^{-1}=M^d
\)
we get
\(
\det(M)=\det(M)^d.
\)
Since \(\det(M)\in 1+\mathfrak m_R\) and \(d-1\in R^\times\), it follows that
\(
\det(M)=1.
\)
Put $N:=M-I_2$, and put
\[
s:=\operatorname{tr}(M)-2=\operatorname{tr}(N)\in\mathfrak m_R.
\]
Using \(\det(M)=\det(I+N)=1\), we have
\(
1=\det(N+I_2)=t^2-\operatorname{tr}(N)t+\det(N)|_{t=-1}=1+\operatorname{tr}(N)+\det(N)
\)
and hence $\det(N)=-s$. Therefore, we have
\(
N^2-sN-sI=0,
\)
or equivalently
\(
N^2=s(N+I).
\)
It follows that
\[
\operatorname{tr}(M^d)-\operatorname{tr}(M)
=
s\bigl((d^2-1)+sW_d(s)\bigr)
\]
for some \(W_d(s)\in R\hoge[s]\). We will explain this $W_d(s)$ shortly. Define polynomials \(P_n(s)\in \mathbb Z[s]\) by
\(
P_0(s)=2,\ \ P_1(s)=s,
\)
and
\(
P_{n+2}(s)=sP_{n+1}(s)+sP_n(s),  \ (n\ge0).
\)
Then, from \(N^2=s(N+I)\), we have
\(
\operatorname{tr}(N^n)=P_n(s)
\)
for all \(n\ge0\). In particular,
\(
P_2(s)=s^2+2s,
\)
and the recurrence gives
\(
P_n(s)\in s^2\mathbb Z[s]\ \ (n\ge3).
\)
Hence
\[
\operatorname{tr}(M^d)
=
\sum_{n\ge0}\binom{d}{n}P_n(s)
=
2+ds+\binom d2(2s+s^2)
+\sum_{n\ge3}\binom dn P_n(s).
\]
Therefore
\[
\operatorname{tr}(M^d)
=
2+d^2s+s^2W_d(s)
\]
for some \(W_d(s)\in \mathbb Z_p\hoge[s]\). Since
\(
\operatorname{tr}(M)=2+s,
\)
we get
\[
\operatorname{tr}(M^d)-\operatorname{tr}(M)
=
(d^2-1)s+s^2W_d(s).
\]
Thus
\[
\operatorname{tr}(M^d)-\operatorname{tr}(M)
=
s\bigl((d^2-1)+sW_d(s)\bigr).
\]
Since \(d\in\mu_{p-1}\) and \(d\neq\pm1\), the element \(d^2-1\) is a unit. Hence
\(
(d^2-1)+sW_d(s)\in R^\times.
\)
Therefore the equality
\(
\operatorname{tr}(M^d)=\operatorname{tr}(M)
\)
implies
\(
s=0.
\)
Thus
\(
\operatorname{tr}(M)=2,\ \det(M)=1.
\)
By Cayley--Hamilton for \(M\), we obtain
\(
(M-I)^2=0.
\) Let $N:=M-I$. Then, $N^2=0$ and hence $M^{\chi(g)}=(I+N)^{\chi(g)}=I+\chi(g)N$. 

If we write $N=\begin{pmatrix}
r& a\\
b& -r\end{pmatrix}$, we obtain
$$\begin{pmatrix}
1+r& \chi_0(g)a\\
\chi_0^{-1}(g)b& 1-r
\end{pmatrix}=\begin{pmatrix}
1+\chi(g)r&\chi(g)a\\
\chi(g)b& 1-\chi(g)r
\end{pmatrix}$$
In general, for each $1\neq x\in \mu_{p-1}$, we know that $x-1\in \Z_p^\times$ because, if not, $x\equiv 1 \mod p\Z_p$. Hence, $r=0$ if $\chi\neq 1$, and $a=b=0$ if $\chi\neq \chi_0,\chi_0^{-1}$.
When $\bar\rho(s_{\chi,i})=\begin{pmatrix}1& 1\\
0& 1\end{pmatrix}$, then we have $\chi=\chi_0$. For such $s_{\chi,i}$, we always denote $s_{\chi,i}=s_{\chi_0,1}$, and denote $M:=\rho_G(s_{\chi_0,1})$.
Then, 
\(
\bar M=
\begin{pmatrix}
1&1\\
0&1
\end{pmatrix}
\in \mathrm{GL}_2(R/\mathfrak m_R), 
\)
and  we may write
\[
M=I+N,\qquad
N=E_{12}+N_0,\qquad N_0\in \mathfrak m_R M_2(R).
\]
Although \(N\) is not contained in \(\mathfrak m_RM_2(R)\), it is still topologically nilpotent, since
\(
\bar N^2=E_{12}^2=0.
\)
Hence \(N^2\in \mathfrak m_RM_2(R)\), and consequently \(N^n\to 0\) \(\mathfrak m_R\)-adically. Therefore, for every \(a\in\mathbb Z_p\), the binomial expression
\[
M^a=(I+N)^a:=\sum_{n\ge0}\binom{a}{n}N^n
\]
is well-defined.

Thus the same argument as in the case \(M\in \ker(\mathrm{GL}_2(R)\to \mathrm{GL}_2(R/\mathfrak m_R))\) applies, with \(N\) replaced by
\(
N=E_{12}+N_0,
\) We obtain that $M=\begin{pmatrix}1&1+X\\
0& 1\end{pmatrix}$ with some $X\in m_R$. Conjugating by an element of $\ker(\mathrm{GL}_2(R)\to \mathrm{GL}_2(\F_p))$, which does not change the deformation class, we may assume that $X=0$. Hence, $\rho_G(s_{\chi_0,1})=\begin{pmatrix}1&1\\
0& 1\end{pmatrix}$.

Finally, for $\chi=\chi_0^{-1}$, $\rho_G(s_{\chi,i})=\begin{pmatrix}1&0\\
Y& 1\end{pmatrix}$ commutes with $\begin{pmatrix}1&1\\
0& 1\end{pmatrix}$, and hence, $Y=0$. For the trivial character component, the above argument shows that
\(\rho_G(s_1)\) is diagonal. Since \(s_1\) and \(s_{\chi_0,1}\) lie in the
abelian group \(X_F\), they commute. We have normalized
\[
\rho_G(s_{\chi_0,1})=
\begin{pmatrix}
1&1\\
0&1
\end{pmatrix}.
\]
A diagonal matrix commuting with this unipotent matrix is scalar. Hence
\[
\rho_G(s_1)=
\begin{pmatrix}
1+V&0\\
0&1+V
\end{pmatrix}.
\]
\\
Step 2. Proof in the case when $\bullet =\pi$.

Note that $G=\mathfrak{X}_{F,\pi}\rtimes (\Delta\times \Gamma)$ when $\bullet=\pi$. Hence, all the relations are given by the relations of $\Lambda$-modules of $\mathfrak{X}_{F,\pi}$.
We have the projective resolution 
$$0\rightarrow \bigoplus_{i=1}^{s_\chi}\Lambda e_{\chi,i}\stackrel{(f_{k,\ell,\chi})_{k,\ell}}{\rightarrow}\bigoplus_{i=1}^{s_\chi}\Lambda e_{\chi,i}\stackrel{\psi}{\rightarrow} \mathfrak{X}_{F,\pi}^{\chi}\rightarrow 0$$
where $\Psi(e_{\chi,i})=s_{\chi,i}$. 
Moreover, $ \rho_G(s_{\chi_0,i}^{\gamma-1})=\begin{pmatrix}1&(\frac{1+S}{1+T}-1)U_i\\
0& 1\end{pmatrix}$.
Hence, for $(f_{k,\ell}):=(f_{k,\ell,\chi_0})$, we obtain
$$f_{1,\ell}(\frac{1+S}{1+T}-1)+\sum_{k=2}^{s_{\chi_0}}U_kf_{k,\ell}(\frac{1+S}{1+T}-1)=0,$$
$\ell=1,\ldots, s_{\chi_0}$ in $R_G$. 
The entries $f_{k,\ell}$ form a presentation matrix for the $\chi_0$-component of $\mathfrak{X}_{F,\pi}$. Hence, after extension of scalars to $\Lambda_W$, the ideal generated by the determinants of $(f_{k,\ell})$ is $(\mu_F(\chi_0))$. Finally, for $s_1$, we know that the trivial character component of $\mathfrak{X}_{F,\pi}$ is isomorphic to $\mathfrak{X}_{K,\pi}$, which is cyclic as a $\Lambda$-module with annihilator $\mu_K$.
Hence, we obtain $(1+V)^{a_\pi}-1=0$ in $R_G$. These are all the relations among $S, T, U_i, V$.\\
Step 3. Proof in the case when $\bullet=p$.

In this case, all the relations among the generators in $G$ are given by $[\gamma_1,\gamma_2]=c\in \mathfrak{X}_{F,p}$, and the $\Lambda$-module relations of $\mathfrak{X}_{F,p}$ by Lemma \ref{lemma:two-var pro-p}. Here, $c\in \mathfrak{X}_{F,p}$ is the extension class of $G$.

However, as proved in Step $1$, $\rho_G(\gamma_1)$ and $\rho_G(\gamma_2)$ are diagonal and therefore commute. Hence, we only need to consider the relations of $\mathfrak{X}_{F,p}$ as a $\Lambda$-module. For the trivial character component, we have
\[
\mathfrak{X}_{K,p}=\Lambda\cdot s_1+\Lambda\cdot c,
\qquad c=[\gamma_1,\gamma_2].
\]
Choose a presentation
\[
0\longrightarrow \Lambda
\xrightarrow{(f,g)}
\Lambda^2
\longrightarrow \mathfrak{X}_{K,p}
\longrightarrow 0,
\]
where the two standard basis vectors of \(\Lambda^2\) map to \(s_1\) and
\(c\), respectively.  Thus, the corresponding multiplicative relation in
the group is
\[
s_1^{\,f(\gamma_1-1,\gamma_2-1)}
\,
c^{\,g(\gamma_1-1,\gamma_2-1)}
=1.
\]
Applying the universal deformation, we have
\[
\rho_G(c)=\rho_G([\gamma_1,\gamma_2])=1,
\]
because \(\rho_G(\gamma_1)\) and \(\rho_G(\gamma_2)\) are diagonal and
therefore commute.  Since
\[
\rho_G(s_1)=
\begin{pmatrix}
1+V&0\\
0&1+V
\end{pmatrix},
\]
the above relation gives
\[
(1+V)^{f(0,0)}-1=0.
\]
We put
\(
a:=f(0,0)\in\mathbb Z_p.
\)
Let $G_\chi$ be a presentation matrix for $\mathfrak{X}_{F,p}^\chi$. By Corollary~\ref{cor:p-iwasawa}, after extension of scalars to $\Lambda_W$, the ideal generated by the determinants of the maximal submatrices of $G_\chi$ contains
\[
\mathfrak{a}\mu_{F,\pi}(\chi)+\mathfrak{a}^{'}\mu_{F,\bar\pi}(\chi)
\]
for some non-zero ideals $\mathfrak{a},\mathfrak{a}^{'}\subset \Lambda_W$.

The remaining arguments are the same as Step 2 by combining Corollary \ref{cor:p-iwasawa},  so we omit the details, but we point out that $\mathfrak{X}_{F,p}^{\chi}$ has $\Lambda$-rank $1$ for all $\chi\in \widehat{\Delta}$. We already proved this in the proof of Proposition \ref{prop:p-iwasawa-2}.
\end{proof}
From this proof, we obtain the following corollary.
\begin{corollary}\label{cor:cyclicity}
We use the terminology of Theorem \ref{Main 2}. Fix a prime $\mathfrak{P}$ of $F$ above $\bar\pi$. 
The ideal $\mathfrak{a}\subset \Lambda_W$ in Theorem \ref{Main 2} is $(1)$ when  either
\(
        F_{\mathfrak P}\not\supset K_{\bar{\pi}}(\zeta_p),
\)
or
\[
        F_{\mathfrak P}\supset K_{\bar{\pi}}(\zeta_p)
        \quad\text{and}\quad
        \chi|_{\Delta_{\mathfrak P}}\neq \omega_{\mathfrak P},
\]
where \(\omega_{\mathfrak P}\) denotes the local mod-\(p\) cyclotomic
character on \(\Delta_{\mathfrak P}:=\Gal(F_\mathfrak{P}/K_{\bar\pi})\). 
\end{corollary}
\begin{proof}
This is immediate from Proposition \ref{prop:p-iwasawa-2}.
\end{proof}
\begin{corollary}\label{cor:free-deform}
Let $\bar\rho:G_K^{\bullet}\rightarrow \mathrm{GL}_2(\F_p)$ , $\bullet\in \{p,\pi\}$ be the reducible Galois representation as above. Let $\rho^{\mathrm{univ}}:G_K^{\bullet}\rightarrow \mathrm{GL}_2(R_{\bar\rho})$ be the universal deformation. For the elements $s_{\chi_0^{-1},1},\ldots, s_{\chi_{0}^{-1},s_{\chi_0^{-1}}}$,   
$$\rho^{\mathrm{univ}}(s_{\chi_0^{-1},i})=\begin{pmatrix}
1&0\\
V_i& 1\end{pmatrix}$$
for some $V_i\in m_{R_{\bar\rho}}$. Furthermore, in the $p$-ramified case, let $F$ be the maximal abelian subextension of $\overline{K}^{\mathrm{ker}(\bar\rho)}$ and assume that the Galois group of the maximal pro-$p$ extension unramified outside $p$ over $F$ is a free pro-$p$ group. Then, the universal deformation ring is $R_{\bar\rho}=\Z_p\hoge[S_1,T_1,S_2,T_2,Y]$. 
\end{corollary}
\begin{proof}
The first half claim is already proved in the proof of Theorem \ref{Main 2}, so we only prove the remaining claim.
Let $G_{F,p}$ be the Galois group of the maximal pro-$p$ extension of $F$ unramified outside $p$ and let $d:= \dim_{\F_p}H^1(G_{F,p}, \F_p)$ be the generator rank of $G_{F,p}$. Moreover, $H^2(G_{F,p},\F_p)=0$ because $G_{F,p}$ is a free pro-$p$ group by the assumption.
By the Euler-characteristic formula, we have 
$0-d+1=-[F:K]$, and hence $d=[F:K]+1$. Now, let $F_\infty$ be the maximal $\Z_p^2$-extension of $F$ unramified outside $p$ which is abelian over $K$, Then, the Galois group $\mathfrak{X}_{F,p}$ of maximal abelian pro-$p$ extension of $F_\infty$ unramified outside $p$ over $F_\infty$ is a free $\Lambda$-module with rank $d-1=[F:K]$. See, Proposition \ref{thm:projective-resolution}. On the other hand, $\mathfrak{X}_{F,p}^{chi}$ has $\Lambda$-rank $1$ for each $\chi\in\widehat{\Delta}$-component (see the proof of Theorem\ref{Main 2}). In particular, we have topological generators $\kakko<\gamma_1,\gamma_2, s_{\chi,1}\ |\ 1\neq \chi\in \widehat{\Delta}>$ by means of Definition \ref{def:Special-gen}. Combining the fact that $G_{F,p}$ is a free pro-$p$ group with the discussion in the proof of Theorem \ref{Main 2} and the first half claim of this corollary, we know that $\rho^{\mathrm{univ}}(s_{\chi,1})=1$ if $\chi\neq \chi_0,\chi_0^{-1}, 1$ and 
$$\rho^{\mathrm{univ}}(s_{\chi_0,1})=\begin{pmatrix}1&1\\
0&1\end{pmatrix},\ \rho^{\mathrm{univ}}(\gamma_i)=\begin{pmatrix}1+S_i& 0\\
0 & 1+T_i\end{pmatrix}\ ,\ \rho^{\mathrm{univ}}(s_{\chi_0^{-1},1})=\begin{pmatrix}1& 0\\
Y& 1\end{pmatrix}.
$$
$G_{F,p}$ is a free pro-$p$ group and hence, no relation appears in the deformation ring, hence we obtain our result.
\end{proof}
\begin{corollary}\label{cor:isom-condition}
Assume that $\mathrm{Cl}_K\otimes \Z_p=0$ and $p$ splits in $K/\Q$.
Let $M/F_\infty$ be the maximal abelian pro-$p$ extension unramified outside $\bullet$ as before for each $\bullet \in \{p,\pi\}$. Let $N:=\Gal(F^\bullet(p)/M)$, and $P_F:=\Gal(M/K)$. Then, the following are equivalent;
\begin{enumerate}
\item The universal deformation $\rho^{\mathrm{univ}}:G_K^{\bullet}\rightarrow \mathrm{GL}_2(R_{\bar\rho})$ is reducible
\item $R_{\bar\rho}\twoheadrightarrow R_G$ is an isomorphism.
\item For the elements $s_{\chi_0^{-1},1},\ldots, s_{\chi_0^{-1}, s_{\chi_0^{-1}}}$, $\rho^{\mathrm{univ}}(s_{\chi_0^{-1},j})=I_2$ for all $j$.
\end{enumerate}
\end{corollary}
\begin{proof}
Step $1$. we prove 
$(1)\Leftrightarrow (2)$.

The implication (2) $\Rightarrow$ (1) is immediate from Theorem \ref{Main 1}. Assume (1) holds. Then, the fixed field $L$ of $\rho^{\mathrm{univ}}$ is contained in the maximal two-step solvable extension $K^{2-\mathrm{step}}$ of $K$ and also contained in the maximal pro-$p$ extension $F^\bullet(p)$ of $F$ unramified outside $\bullet$. Then $L\subset K^{2-\mathrm{step}}$ and pro-$p$ over $F$. Let $F_\infty:=K^{\mathrm{ab}}\cap L$. Then, $L/F_\infty$ is abelian because $K^{2-\mathrm{step}}/K^{\mathrm{ab}}$ is abelian. Note that $F_\infty/F$ is an abelian pro-$p$ extension and $F_\infty/K$ is abelian. Let $K_\infty(p)$ be the maximal abelian pro-$p$ extension of $K$ unramified outside $\bullet$. From this we obtain $F_\infty\subset K_\infty(p)\cdot F$, and by our assumption that $\mathrm{Cl}_K\otimes \Z_p=0$, we obtain $K_\infty(p)/K$ is a multiple $\Z_p$-extension. Hence, the representation $\rho^{\mathrm{univ}}$ factors through $G^\bullet$.

Step $2$. we prove $(3)\Leftrightarrow (1)$. 

$(3)$ implies $(1)$ is obvious, so we will show $(1)$ implies $(3)$. Assume that there exists some $j$ such that $\rho^{\mathrm{univ}}(s_{\chi_0^{-1},j})\neq I_2$. 
Then, there exist elements $a,b$ such that $$\rho^{\mathrm{univ}}(a)=\begin{pmatrix}1& 1\\
0& 1\end{pmatrix}, \rho^{\mathrm{univ}}(b)=\begin{pmatrix} 
1& 0\\
z& 1\end{pmatrix},$$
with some $z\in m_{R_{\bar\rho}}-\{0\}$. In the following, we simply note $\rho=\rho^{\mathrm{univ}}$.

Suppose that \(L\subset V:=R_{\bar\rho}e_1\oplus R_{\bar\rho}e_2\) is a
\(G\)-stable rank one direct summand. Since \(L\) is a direct summand,
\(\overline L:=L\otimes_{R_{\bar\rho}} \F_p\) is a one-dimensional subspace of
\(V\otimes_{R_{\bar\rho}} \F_p\). It is stable under
\(\overline{\rho(a)}=\begin{pmatrix}1&1\\0&1\end{pmatrix}\), whose unique
invariant line over \(\F_p\) is \(\F_p\overline e_1\). Hence \(L\) is generated by a
vector of the form \(e_1+t e_2\) with \(t\in\mathfrak m_{R_{\bar\rho}}\).

Since \(L\) is stable under \(b\), we have
\(\rho(b)(e_1+t e_2)\in L\). But
\(\rho(b)(e_1+t e_2)=e_1+(z+t)e_2\). Thus there exists
\(\lambda\in R_{\bar\rho}\) such that
\(e_1+(z+t)e_2=\lambda(e_1+t e_2)\). Comparing the first coordinate gives
\(\lambda=1\), and comparing the second coordinate gives \(z+t=t\).
Therefore \(z=0\), contradicting the choice of \(z\). Hence no such \(L\)
exists.
\end{proof}
In the next section, we treat a Galois representation whose ramification set contains the primes above $pd_K$. For this purpose, we enlarge the ramification set in Theorem \ref{Main 1}.
\begin{corollary}
Let $S$ be a finite set of finite places of $K$. Assume that every $\mathfrak{q}\in S$ satisfies $\mathrm{N}\mathfrak{q}\not\equiv 1 \mod p$. Then Theorem \ref{Main 2} remains valid for the reducible Galois representation $\bar\rho:G_K^S\rightarrow \mathrm{GL}_2(\F_p)$ with noncommutative image and nonquadratic character $\chi_0\neq 1$.
\end{corollary}
\begin{proof}
Put
\(
S^{(p)}:=S\cap S_p.
\)
We prove that allowing ramification at the primes in
\(
S\setminus S^{(p)}
\)
does not change the maximal abelian pro-$p$ extension of $K_\infty$.

Let $\mathfrak q\in S\setminus S^{(p)}$, and let $w$ be a prime of $K_\infty$
above $\mathfrak q$. Since $\mathfrak q\nmid p$, the extension
\(
K_{\infty,w}/K_{\mathfrak q}
\)
is unramified. Hence, for each finite layer $K_n/K$ and each prime $w_n$ of
$K_n$ below $w$, the local extension
\(
K_{n,w_n}/K_{\mathfrak q}
\)
is unramified of $p$-power degree. Therefore its residue field has cardinality
\(
Nw_n=(N\mathfrak q)^{p^a}
\)
for some $a\geq 0$.

By assumption,
\(
N\mathfrak q\not\equiv 1 \mod p.
\)
Since $p$ is odd, this implies
\(
(N\mathfrak q)^{p^a}\not\equiv 1\mod p.
\)
Thus
\(
\mu_p(K_{n,w_n})=1
\)
for every $n$. Passing to the union, we obtain
\(
\mu_{p^\infty}(K_{\infty,w})=1.
\)

Now we use local class field theory. Let $E$ be a finite extension of
$K_{\mathfrak q}$ contained in $K_{\infty,w}$. Since the residue characteristic
of $E$ is different from $p$, the pro-$p$ completion of $E^\times$ has the form
\(
\widehat{E^\times}
\simeq
\Z_p\oplus \mu_{p^\infty}(E).
\)
The $\Z_p$-factor corresponds to the unramified $\Z_p$-extension of $E$.
Therefore, if
\(
\mu_{p^\infty}(E)=1,
\)
then $E$ has no nontrivial ramified abelian pro-$p$ extension.

Applying this to all finite layers $K_{n,w_n}$ and passing to the direct limit,
we see that $K_{\infty,w}$ has no nontrivial ramified abelian pro-$p$ extension.
Hence every abelian pro-$p$ extension of $K_\infty$ which is allowed to ramify
at $w$ is in fact unramified at $w$.

Since this holds for every prime $w$ of $K_\infty$ above every
\(
\mathfrak q\in S\setminus S^{(p)},
\)
the maximal abelian pro-$p$ extension of $K_\infty$ unramified outside $S$ is
already unramified outside $S^{(p)}$. Hence
\(
M_S(K_\infty)\subset M_{S^{(p)}}(K_\infty).
\)
The reverse inclusion is immediate from
\(
S^{(p)}\subset S.
\)
Therefore
\(
M_S(K_\infty)=M_{S^{(p)}}(K_\infty).
\)
Taking Galois groups over $K_\infty$, we obtain
\(
\mathfrak X_{K,S}
=
\Gal(M_S(K_\infty)/K_\infty)
\simeq
\Gal(M_{S^{(p)}}(K_\infty)/K_\infty)
=
\mathfrak X_{K, S^{(p)}}.
\)
In particular, when $S^{(p)}=S_p$, this gives
\(
\mathfrak X_{K,S}\simeq \mathfrak X_{K,S_p},
\)
and when $S^{(p)}=S_\pi$, it gives
\(
\mathfrak{X}_{K,S}\simeq \mathfrak{X}_{K, S_\pi}.
\)
\end{proof}
\section{Universal Galois representations attached to Bianchi cusp forms}
In this section, we study the Galois representation of Hecke-eigen Bianchi cusp forms with cyclotomic central character, studied by Harris--Soudry--Taylor\cite{HarrisSoudryTaylor1993} and Berger--Klosin\cite{BergerKlosin}. In the following, we fix two field embeddings, $\iota_\infty:\overline{\Q}\hookrightarrow \C$, $\iota_p:\overline{\Q}\hookrightarrow \C_p$.

Let \(K\) be an imaginary quadratic field, and let
\(\Pi=\otimes_v'\Pi_v\) be a cuspidal automorphic representation of
\(\mathrm{GL}_2(\mathbb A_K)\) attached to a Bianchi newform \(f\) of level $K_f\subset \mathrm{GL}_2(\A_{K,\mathrm{fin}})$
with the central character
\(\psi\). For the definition of this automorphic representation, see \cite[Section 3.3]{Bump}. Let \(E_f\) be the Hecke field of \(f\), and let \(\lambda\) be a
finite prime of \(E_f\) above a rational prime \(p\). Assume that the following two conditions hold;
\begin{enumerate}
\item $\Pi_\infty$ is a principal-series $(\mathfrak{g},\mathrm{SU}(2))$ module obtained by Borel induction of 
$$\begin{pmatrix}
t_1& *\\
0& t_2\end{pmatrix}\mapsto \frac{t_1}{|t_1|}\frac{|t_2|}{t_2}.$$
Here, $\mathfrak{g}$ is the Lie algebra of the real Lie group $\mathrm{GL}_2(\R\otimes K)$.
\item $\psi$ is the cyclotomic character.
\end{enumerate}
Let $S$ be a set of finite places of $K$ containing the places above $pd_K$ and the finite primes \(\mathfrak q\) for which
\(\Pi_{\mathfrak q}\) is not spherical, equivalently
\(
\Pi_{\mathfrak q}^{\mathrm{GL}_2(\mathcal O_{K_{\mathfrak q}})}= 0.
\)

Then there is a
continuous absolutely irreducible Galois representation
\[
\rho_{f,\lambda}:G_K^{S}\longrightarrow \mathrm{GL}_2(E_{f,\lambda})
\]
which is unramified outside \(S\). For every finite prime
\(\mathfrak q\not\in S\), it satisfies
\[
\operatorname{tr}\rho_{f,\lambda}(\mathrm{Frob}_{\mathfrak q})
=a_{\mathfrak q}(f),
\qquad
\det\rho_{f,\lambda}(\mathrm{Frob}_{\mathfrak q})
=\psi(\mathfrak q)\mathrm{N}\mathfrak q,
\]
Here \(a_{\mathfrak q}(f)\) denotes the Hecke eigenvalue of \(f\) at
\(\mathfrak q\) and $\mathrm{N}\mathfrak{q}:=|\O_K/\mathfrak{q}|$.

Let \(\mathcal O_\lambda\) be the valuation ring of \(E_{f,\lambda}\), and
let \(k_\lambda\) be its residue field. After choosing a stable lattice, we
obtain a residual representation
\[
\bar\rho_{f,\lambda}:G_K^S\longrightarrow \mathrm{GL}_2(k_\lambda).
\]
The semisimplification of \(\bar\rho_{f,\lambda}\) is independent of the
choice of the lattice. We say that \(\lambda\) is an Eisenstein prime for
\(f\) if \(\bar\rho_{f,\lambda}\) is reducible. In this case, one has
\[
\bar\rho_{f,\lambda}^{\mathrm{ss}}
\simeq \bar\chi_1\oplus \bar\chi_2
\]
for characters \(\bar\chi_i:G_K^S\to k_\lambda^\times\). Equivalently, for
almost all primes \(\mathfrak q\) of \(K\),
\[
a_{\mathfrak q}(f)\equiv
\bar\chi_1(\mathrm{Frob}_{\mathfrak q})
+
\bar\chi_2(\mathrm{Frob}_{\mathfrak q})
\mod{\lambda}.
\]
Thus, an Eisenstein congruence gives rise to a reducible residual Galois
representation. By Ribet's lattice lemma, after choosing a
suitable lattice and its basis,  we may write
\[
\bar\rho_{f,\lambda}=
\begin{pmatrix}
\bar\chi_1 & *\\
0 & \bar\chi_2
\end{pmatrix},
\qquad *\neq 0 .
\]
Berger--Klosin\cite{BergerKlosin} proved the following modularity theorem for mod $p$ reducible representations.
\begin{theorem}(Berger--Klosin\cite{BergerKlosin}, Theorem 5.8, Theorem 5.17)\label{Theorem modularity}
Fix a $p$-adic field $E$. For an irreducible representation $\rho:G_K^{S}\rightarrow \mathrm{GL}_2(E)$, assume that $\bar\rho^{ss}\simeq \chi_1\oplus \chi_2$ and put $\chi_0:=\chi_1/\chi_2$. Under the following conditions, $\rho$ is modular.
\begin{enumerate}
\item There exist Hecke characters $\phi_1,\phi_2:\A_K^\times/K^\times\rightarrow \C^\times$ of $\infty$-type $\phi_1(z_\infty)=z_\infty, \phi_2(z_\infty)=z_\infty^{-1}$, such that $\chi_1=\overline{\phi_{1,p}\epsilon_{cyc}}, \chi_2=\overline{\phi_{2,p}}$. Here, $\phi_{i,p}:G_K^S\rightarrow E^\times$ is the $p$-adic Galois representation attached to $\phi_i$.
\item $p\nmid\ |(\O_K/\mathrm{cond}(\phi_1))^\times |$,
\item $\frac{L(0, \phi_1/\phi_2)}{\Omega^2}\in E$ and $\mathrm{val}_E(\frac{L(0, \phi_1/\phi_2)}{\Omega^2})=1$, where $\Omega\in \C^\times$ is the complex period.
\item $\chi_0$ and $\chi_0^{-1}$ are ramified at $\pi$, and for all $q\in S$, either $\chi_0$ is ramified at $q$ or $\chi_0(\mathrm{Frob}_q)^{-1}\neq \mathrm{N}(q)$,
\item $\chi_0(cgc)=\chi_0(g)^{-1}$ for all $g\in G_K^S$, $c$: complex conjugation.
\item $\mathrm{Cl}_K\otimes \Z_p=\mathrm{Cl}_{K(p)}\otimes \Z_p=0$.
\item $p\nmid[K(p):K]$.
\item $S$ contains the primes above $pd_K\mathrm{N}(\mathrm{cond}(\phi_1))$, 
\item $\det(\rho)=\phi_1\phi_2\epsilon_{\mathrm{cyc}}$
\item For all $\mathfrak{q}\in S$, $\mathrm{N}\mathfrak{q}\not\equiv 1 \mod p$.
\end{enumerate}
Furthermore, for a given reducible representation $\bar\rho$ satisfying the above, there exists an irreducible deformation of $\bar\rho$.
\end{theorem}
A numerical example of such a representation $\bar\rho:G_K^S\rightarrow \mathrm{GL}_2(k_E)$ is given by Berger--Klosin. 
\begin{prop}\label{prop:autom}
Let $K=\Q(\sqrt{-51})$ and $p=5$. Then $(p)=(5, \omega-2)(5, \omega+1)$ where $\omega:=\frac{1+\sqrt{-51}}{2}$. Put $\mathfrak{p}:=(5, \omega-2)$. Then, there exists a modular Galois representation $\rho:G_K^{\{5,3,17\}}\rightarrow \mathrm{GL}_2(\O_{K,\mathfrak{p}})$ such that 
$$\bar\rho=\begin{pmatrix}
\chi_1 & *\\
0& \chi_2\end{pmatrix},\  *\neq 0$$
and $\chi_0:=\chi_1/\chi_2\neq 1$ is not a quadratic character. We denote by $R_{\bar\rho}$ the universal deformation ring of $\bar\rho$. Let $(3)=\mathfrak{p}_3^2$ and let $F\subset K(5\mathfrak{p}_3)$ be the fixed field of $\chi_0 :G_K^{\{5,3,17\}}\rightarrow \F_5^\times$. Then $R_{\bar\rho}\simeq \Z_5\hoge[S_1,T_1,S_2,T_2,U_2,Y_1,Y_2]/I$ where $I$ is an ideal contained in $(Y_1,Y_2, g_1((1+S_1)/(1+T_1)-1, (1+S_2)/(1+T_2)-1)+U_2g_2((1+S_1)/(1+T_1)-1, (1+S_2)/(1+T_2)-1))$ with some power series $g_1,g_2\in \Z_p\hoge[S,T]$.
\end{prop}
\begin{proof}
First, let $\phi:\A_K^\times/K^\times \rightarrow \C^\times$ be an unramified Hecke character such that $\phi(z_\infty)=z_\infty^{2}$, $\phi|_{\mathrm{fin}}=1$. Note that $\mathrm{Cl}_K=2$ and hence, there are two unramified Hecke characters of infinity type $z_\infty^2$. Berger--Klosin showed that $L(0,\phi)/\Omega^2\in \Z_5$ and its valuation satisfies $\mathrm{val}_5(L(0,\phi)/\Omega^2)=1$, using Magma.

Step $1$. construction of $\phi_1$, $\phi_2$. 

We write
\(
3=\mathfrak p_3^2,\ 
\mathfrak p_3=(3,\omega+1),
\)
and
\(
5=\mathfrak p\bar{\mathfrak p},
\ 
\mathfrak p=(5,\omega-2),\ 
\bar{\mathfrak p}=(5,\omega+1).
\)
Put
\[
\mathfrak m=\mathfrak p_3(5)=\mathfrak p_3\mathfrak p\bar{\mathfrak p}.
\]

Let $I_{\mathfrak m}$ be the group of fractional ideals of $K$ prime to
$\mathfrak m$, and let $P_{\mathfrak m,1}$ be the subgroup generated by principal ideals
$(a)$ with
\(
a\equiv 1 \mod{\mathfrak m}.
\)
We use the following convention for algebraic Hecke characters: a character
$\psi:I_{\mathfrak m}\to \overline{\mathbb Q}^{\times}$ has infinity type $(r,s)$ if
\(
\psi((a))=a^{-r}\bar a^{-s}
\)
for all $(a)\in P_{\mathfrak m,1}$.

Let
\(
\varepsilon_3:I_{\mathfrak m}\longrightarrow \{\pm1\}
\)
be the quadratic ray class character of conductor $\mathfrak P_3$. Equivalently,
$\varepsilon_3$ factors through
\(
(\mathcal O_K/\mathfrak P_3)^\times\simeq \mathbb F_3^\times
\)
on principal ideals prime to $\mathfrak P_3$.

Let also
\(
\eta:I_{\mathfrak m}\longrightarrow \overline{\mathbb Q}^{\times}
\)
be a finite order ray class character of conductor dividing $(5)$, passing through
$$\eta:(\O_K/5\O_K)^\times \simeq (\O_K/\lambda)^\times \times (\O_K/\bar\lambda)^\times \simeq \F_5^\times  \times \F_5^\times\rightarrow \F_5^\times;(a,b)\mapsto \frac{a}{b}$$.

Define two algebraic Hecke characters
\[
\phi_1,\phi_2:I_{\mathfrak m}\longrightarrow \overline{\mathbb Q}^{\times}
\]
by prescribing their values on principal ideals prime to $\mathfrak m$ as follows:
\(
\phi_1((a))
=
a^{-1}\varepsilon_3(a)\eta(a),
\)
and
\(
\phi_2((a))
=
a\,\varepsilon_3(a)\eta(a).
\)
Here $\varepsilon_3(a)$ and $\eta(a)$ denote the finite parts evaluated on the residue class of $a$ modulo their conductors.

Thus $\phi_1$ has infinity type $(1,0)$ and $\phi_2$ has infinity type $(-1,0)$ in the above convention. Their quotient is
\(
\phi:=\frac{\phi_1}{\phi_2}.
\)
On principal ideals one has
\(
\phi((a))
=
\frac{\phi_1((a))}{\phi_2((a))}
=
a^{-2}.
\)
Hence $\phi$ has infinity type $(2,0)$ in the same convention. In the automorphic convention in which infinity type is written by the archimedean function $z^2$, this is the Hecke character of infinity type $z^2$.

We now pass to $5$-adic Galois characters. Fix the $5$-adic embedding corresponding to
\(
\mathfrak p=(5,\omega-2).
\)
Let
\(
\phi_{i,p}:G_K^{\{3,5,17\}}\longrightarrow \Z_5^\times
\)
be the $5$-adic Galois character associated to $\phi_i$ by class field theory. Reducing modulo $5$, we obtain
\(
\bar\phi_{i,p}:G_K\longrightarrow \F_5^\times.
\)

Let
\(
\bar\varepsilon:G_K\longrightarrow \F_5^\times
\)
be the mod $5$ cyclotomic character. The residual Eisenstein representation appearing in the Berger--Klosin construction has semisimplification
\[
\bar\rho^{\mathrm{ss}}
\simeq
\bar\phi_{1,p}\bar\varepsilon
\oplus
\bar\phi_{2,p}.
\]
Equivalently,
\[
\bar\rho^{\mathrm{ss}}
\simeq
\begin{pmatrix}
\bar\phi_{1,p}\bar\varepsilon & 0\\
0 & \bar\phi_{2,p}
\end{pmatrix}.
\]

Twisting by $\bar\phi_{2,p}^{-1}$ gives
\[
\bar\rho^{\mathrm{ss}}\otimes \bar\phi_{2,p}^{-1}
\simeq
\bar\phi_p\bar\varepsilon\oplus 1,
\]
where
\[
\bar\phi_p
=
\bar\phi_{1,p}\bar\phi_{2,p}^{-1}.
\]
Thus the essential character is
\(
\bar\phi_p\bar\varepsilon.
\)

Since
\(
5=\mathfrak p\bar{\mathfrak p},
\)
we have
\(
\bar\varepsilon=\bar\chi_{\mathfrak p}\bar\chi_{\bar{\mathfrak p}},
\)
where
\(
\bar\chi_{\mathfrak p},\bar\chi_{\bar{\mathfrak p}}:G_K\to \mathbb F_5^\times
\)
are the two mod $5$ local characters corresponding to the two primes above $5$.

With respect to the above choice of embedding, the reductions of $\phi_1$ and $\phi_2$ are of the form
\(
\bar\phi_{1,p}
=
\bar\varepsilon_3\,\bar\eta\,\bar\chi_{\mathfrak p},
\)
and
\(
\bar\phi_{2,p}
=
\bar\varepsilon_3\,\bar\eta\,\bar\chi_{\mathfrak p}^{-1}.
\)
Therefore, the two diagonal characters are
\[
\chi_1:=\bar\phi_{1,p}\bar\varepsilon
=
\bar\varepsilon_3\,\bar\eta\,
\bar\chi_{\mathfrak p}^{2}\bar\chi_{\bar{\mathfrak p}},
\]
and
\[
\chi_2:=\bar\phi_{2,p}
=
\bar\varepsilon_3\,\bar\eta\,
\bar\chi_{\mathfrak p}^{-1}.
\]
Their quotient is
\(
\chi_0:=\chi_1/\chi_2=\bar\phi_p\bar\varepsilon
=
\bar\chi_{\mathfrak p}^{3}\bar\chi_{\bar{\mathfrak p}},
\)
which has an image of order $4$. Indeed, after restriction to the inertia group at $\mathfrak p$, the character $\bar\chi_{\bar{\mathfrak p}}$ is trivial, while $\bar\chi_{\mathfrak p}$ has image $(\mathcal O_K/\mathfrak p)^\times\simeq\mathbb F_5^\times$. Hence this quotient has image of order $4$, and in particular $\chi_0$ is not quadratic. The fixed field $F$ of the twisted residual character has degree
\(
[F:K]=4.
\) 

Step $2$. Calculating the universal deformation

The universal deformation ring is unchanged under twisting. Hence, it is  sufficient to consider the universal deformation of $\bar\rho=\begin{pmatrix}1& *\\
0& \chi_0 \end{pmatrix}$.
We also note that the Galois group $G_{F,5}$ of the maximal pro-$5$
extension of $F$ unramified outside the primes above $5$ has seven
topological generators, namely
$\gamma_1,\gamma_2,s_{\chi_0^{-1},1}, s_{\chi_0^{-1},2},s_{\chi_0,1},s_{\chi_0,2}, s_{\chi_0^2,1}$ by means of Definition \ref{def:Special-gen}.
In fact, $\mathrm{Cl}_F\otimes \Z_5=0$ by MAGMA computation, and $\mu_5(F_v)=\mu_5\simeq \Z/5\Z$ for the two primes $v|5$, and hence the Galois group of the pro-$5$ extension of $F$ unramified outside $5$ is isomorphic to $(\prod_{v|5}\O_{F,v}^\times/\overline{\O_F^\times}) \widehat{\otimes}\Z_5\simeq \Z_5^5\times (\Z/5\Z)^{\oplus 2}$. Therefore, $G_{F,5}$ has the $7$ topological generators and hence we can apply Theorem \ref{Main 2} and Corollary \ref{cor:free-deform}.
\end{proof}
\bibliographystyle{alpha}
\bibliography{references}
\end{document}